\documentclass{jams-l}

\usepackage{amssymb}

\usepackage{graphicx}

\usepackage{amsthm}
\usepackage{amssymb}
\usepackage{amsmath}
\usepackage{amsfonts,color}
\usepackage[dvips]{epsfig}
\newcommand{\baseRing}[1]{\ensuremath{\mathbb{#1}}}
\newcommand{\Z}{\baseRing{Z}}
\newcommand{\C}{\baseRing{C}} 
\newcommand{\N}{\baseRing{N}}
\newcommand{\R}{\baseRing{R}}
\newcommand{\Q}{\baseRing{Q}}
\newcommand{\T}{\baseRing{T}}
\newcommand{\D}{\Delta}

\def\P{\baseRing{P}}

\newcommand{\BB}{\mathcal{B}}

\newcommand{\LL}{\mathcal{L}}

\newtheorem{theorem}{Theorem}[section]
\newtheorem{lemma}[theorem]{Lemma}

\newtheorem{corollary}[theorem]{Corollary}
\newtheorem{proposition}[theorem]{Proposition}

\newtheorem{conjecture}[theorem]{Conjecture}

\theoremstyle{definition}

\newtheorem{example}[theorem]{Example}

\theoremstyle{remark}
\newtheorem{remark}[theorem]{Remark}

\numberwithin{equation}{section}

\begin{document}
\title{Tropical Discriminants}


\author{Alicia Dickenstein}
\address{Dto.\ de Matem\'atica, FCEN, Universidad de Buenos Aires,
(1428) B.\ Aires, Argentina}
\email{alidick@dm.uba.ar}
\thanks{A.~Dickenstein was partially supported  by
UBACYT X042, CONICET PIP 5617 and ANPCYT 17-20569, Argentina.}

\author{Eva Maria Feichtner}
\address{Department of Mathematics, ETH Z\"urich, 8092 Z\"urich, Switzerland}
\curraddr{Dept.\ of Mathematics, University of Stuttgart, 70569 Stuttgart, Germany}
\email{feichtne@igt.uni-stuttgart.de}
\thanks{E.M. Feichtner was supported by a Research Professorship of the Swiss 
National Science Foundation, PP002--106403/1.}

\author{Bernd Sturmfels}
\address{Department of Mathematics, University~of California, 
Berkeley CA 94720, USA}
\email{bernd@math.berkeley.edu}
\thanks{B.~Sturmfels was partially supported by
the U.S.~National Science Foundation, DMS-0456960.}


\date{December 18, 2006}

\dedicatory{Dedicated to the memory of Pilar Pis\'on Casares}


\begin{abstract}
Tropical geometry is used to develop a new approach to the theory of
  discriminants and resultants in the sense of Gel'fand, Kapranov and
  Zelevinsky.  The tropical $A$-discriminant is the
  tropicalization of the dual variety of the projective toric variety
  given by an integer matrix~$A$. This tropical algebraic variety is shown to coincide with the
  Minkowski sum of the row space of $A$ and the tropicalization of
  the kernel of~$A$.  This leads to an explicit positive formula for all 
  the extreme monomials of any $A$-discriminant.
\end{abstract}

\maketitle

\section{Introduction}

Let $A$ be an integer $d \times n$-matrix such that $(1,1, \ldots, 1)$
is in the row span of~$A$. This defines a {\em projective toric
  variety} $X_A$ in $\C \P^{n-1}$. Its {\em dual variety} $X_A^*$ is
the closure in the projective space dual to $\C \P^{n-1}$ of the set
of hyperplanes that are tangent to $X_A$ at a regular point.  The
toric variety $X_A$ is called {\em non-defective} if its dual variety
$X_A^*$ has codimension one. In this case, the {\em $A$-discriminant\/}
is the irreducible homogeneous polynomial $\Delta_A$ which defines the
hypersurface $X_A^*$.  Alternatively, the dual variety $X_A^*$ 
can be thought of as the set of singular hypersurfaces in $(\C^*)^d$
with Newton polytope prescribed by the matrix~$A$.
The study of these objects is an active area of
research in computational algebraic geometry, with the fundamental
reference being the monograph by Gel'fand, Kapranov and
Zelevinsky~\cite{GKZ}.

Our main object of interest in this paper is the {\em tropical
  $A$-discriminant\/} $\tau(X_A^*)$. This is the polyhedral fan in
$\R^n$ which is obtained by tropicalizing $X_A^*$. While it is
generally difficult to compute the dual variety $X_A^*$ from $A$, we
show that its tropicalization $\tau(X_A^*)$ can be computed much more
easily.  In Theorem~\ref{thm:tropical}, we derive an explicit
 description of the tropical $A$-discriminant
$\tau(X_A^*)$, and in Theorem~\ref{thm:initial}  we present an explicit
combinatorial formula for the extreme monomials of $\Delta_A$.

Without loss of generality, we assume that 
 the columns of the matrix $A$ span the integer lattice $\Z^d$, and
that the point configuration given by the columns of $A$ is not a
pyramid. These hypotheses ensure that the toric variety $X_A$ has
dimension $d\,{-}\,1$ and that the dual variety $X_A^*$ is not
contained in any coordinate hyperplane.

A key player in this paper is the tropicalization of the kernel of
$A$. As shown in \cite{AK} and \cite{FS}, this tropical linear space
is subdivided both by the {\em Bergman fan\/} of the matroid dual
to~$A$, i.e., the {\em co-Bergman fan\/} $\mathcal{B}^*(A)$ of~$A$,
and by the {\em nested set fans\/} of its lattice of flats~$\LL(A)$.
In other words, let ${\mathcal L}(A)$ denote the geometric lattice
whose elements are the sets of zero-entries of the vectors in ${\rm
kernel}(A)$, ordered by inclusion. We write ${\mathcal C}(A)$ for the
set of proper maximal chains in ${\mathcal L}(A)$ and represent these
chains as $(n\!-\!d\!-\!1)$-element subsets $\sigma =
\{\sigma_1, \dots, \sigma_{n-d-1}\}$ of $\{0,1\}^n$.

We obtain the following descriptions of the tropicalization of the
kernel of $A$:
\begin{equation}
\label{eq:coBergmanFan}
 \tau( {\rm kernel}(A)) \quad = \quad
{\rm support}({\mathcal B}^*(A)) \quad = \quad
\bigcup_{\sigma \in {\mathcal C}(A)}\R_{\geq 0}\, \sigma \, .
\end{equation}
The union on the right hand side in fact is the
finest in the  hierarchy of unimodular simplicial fan 
structures provided by the {\em nested set fans}
\cite{DP, FK, FM, FS, FY}.
The tropical linear space (\ref{eq:coBergmanFan}) is a subset of $\R^n$.
We obtain the tropical $A$-discriminant by adding 
this tropical linear space to the (classical) row space of the
$d \times n$-matrix~$A$:

\begin{theorem}\label{thm:tropical}
For any  $d{\times}n$-matrix $A$ as above, 
the tropical $A$-discriminant $\tau(X_A^*)$ equals
the Minkowski sum of the co-Bergman fan ${\mathcal B}^*(A)$ 
and the row space of~$A$.
\end{theorem}

Theorem \ref{thm:tropical} is the tropical analogue of Kapranov's {\em Horn uniformization\/}~\cite{K}.
By definition, the tropical discriminant $\tau(X_A^*)$ inherits the 
structure of a fan from the Gr\"obner fan of the ideal of $X_A^*$ and, in the
non-defective case, also from the secondary fan of~$A$~\cite{GKZ}. 
In general,
neither of these two fan structures refines the other, as we shall see
in Examples~\ref{ex:K_4} and \ref{ex:BJSST}.

The tropicalization of a variety retains a lot
of information about  the geometry of the original variety
\cite{M,Spe,SS, S, Tev}. 

In Theorem~\ref{thm:initial} below, our tropical approach leads to a
formula for the extreme monomials of the $A$-discriminant~$\Delta_A$, 
and, a fortiori, for the degree of the dual variety $X_A^*$. An
alternating product formula for the extreme monomials of $\Delta_A$
was given in \cite[\S 11.3.C]{GKZ} under the restrictive
assumption that $X_A$ is smooth.  Our formula (\ref{eq:degree1}) is
positive, it is valid for any toric variety $X_A$ regardless of
smoothness, and its proof is self-contained.

\begin{theorem}\label{thm:initial}
  If $X_A$ is non-defective and $w$ a generic vector in $\R^n$ then the
  exponent of $x_i$ in the initial monomial ${\rm in}_w(\Delta_A)$ of
  the $A$-discriminant $\Delta_A$ equals
\begin{equation}\label{eq:degree1}
\sum_{\sigma \in {\mathcal C}_{i,w} } |\, 
      \det ( A^t,\sigma_1, \dots,  \sigma_{n-d-1}, e_i )  \,| \,  ,
\end{equation}
where $\mathcal C_{i,w}$ is the subset of ${\mathcal C}(A)$ consisting of
all chains such that the row space of the matrix $A$ has non-empty
intersection with the cone~$\R_{> 0} \bigl\{\sigma_1, \dots,
\sigma_{n-d-1}, -e_i, -w\bigr\}$.
\end{theorem}

\noindent
Here, the $A$-discriminant $\Delta_A$ is written as a homogeneous
polynomial in the variables $x_1,\ldots,x_n$, and ${\rm
  in}_w(\Delta_A)$ is the $w$-lowest monomial $x_1^{u_1} \cdots
x_n^{u_n}$ which appears in the expansion of $\Delta_A$ in
characteristic zero. Theorem \ref{thm:initial} generalizes to the
defective case, when we take $\Delta_A$ as the {\em Chow form} of the
dual variety $X_A^*$. This is stated in Theorem
\ref{thm:initial2}. Aiming for maximal efficiency in evaluating
(\ref{eq:degree1}) with a computer, we can replace
$\,\mathcal{C}_{i,w}\,$ with the corresponding maximal nested sets of
the geometric lattice $\mathcal L(A)$, or with the corresponding
maximal cones in the Bergman fan $\mathcal{B}^*(A)$. Our {\tt maple}
implementation of the formula (\ref{eq:degree1}) is discussed in
Section~\ref{sec:subdivisions}.

\smallskip This paper is organized as follows. In
Section~\ref{sec:tropical}, we review the construction of the
tropicalization $\tau(Y)$ of a projective variety $Y$, and we show how
the algebraic cycle underlying any initial monomial ideal of $Y$ can
be read off from $\tau(Y)$. In Section~\ref{sec:monomials}, we discuss
general varieties which are parametrized by a linear map followed by a
monomial map. Theorem \ref{thm:tUV} gives a combinatorial description
of the tropicalization of the image of such a map.  The dual variety
$X_A^*$ of any toric variety $X_A$ admits such a parametrization. 
This is derived in Section~\ref{sec:Horn}, and it is used to prove 
Theorem \ref{thm:tropical} and Theorem
\ref{thm:initial} in the general form of Theorem \ref{thm:initial2}.
We also obtain in Corollary~\ref{cor:dimension} 
a characterization of the dimension of $X_A^*$.
In Section~\ref{sec:subdivisions} we discuss computational issues,  and
we examine the connection between the tropical discriminant of $A$ and
regular polyhedral subdivisions. In particular, we consider the problem 
of characterizing $\Delta$-equivalence
of regular triangulations in  combinatorial terms. 
Finally, Section~\ref{sec:resultants} is devoted to the case when $A$
is an {\em essential Cayley configuration}. The corresponding dual
varieties $X_A^*$ are resultant varieties, and we compute their
degrees and initial cycles in terms of mixed subdivisions.

\bigskip

\noindent {\bf Acknowledgement}: We thank the Forschungsinstitut 
f\"ur Mathematik at ETH Z\"urich for hosting Alicia Dickenstein
and Bernd Sturmfels in the summer of 2005.
We are grateful to Jenia Tevelev and Josephine Yu for comments
on this paper.


\section{Tropical varieties and their initial cycles}
\label{sec:tropical}

{\em Tropical algebraic geometry} is algebraic geometry over the
tropical semi-ring $(\R\,{\cup}\,\{\infty\},\oplus,\odot )$ with arithmetic
operations $x\oplus y := {\rm min} \{x,y\}$ and $x\odot  y := x{+}y$.
It transfers the objects of classical algebraic geometry into the
combinatorial context of polyhedral geometry. 
Fundamental  references include
\cite{EKL, M, RST, Spe, SS, Tev}.

Tropicalization is an operation that turns complex projective
varieties into polyhedral fans. Let $Y\,\subset\,\C\P^{n-1}$ be an
irreducible projective variety of dimension $r-1$ and $I_Y\subset
\C[x_1,\ldots,x_n]$ its homogeneous prime ideal. For $w\in \R^n$, we denote by
${\rm in}_{w}(I_Y)$ the initial ideal generated by all initial forms
${\rm in}_{w}(f)$, for $f\,{\in}\,I_Y$, where ${\rm in}_{w}(f)$ is the
subsum of all terms $c_ux^u\,{=}\,c_u\,x_1^{u_1}\cdots 
x_n^{u_n}$ in~$f$ which have {\em lowest\/} weight
$w\,{\cdot}\,u\,{=}\,\sum_{i=1}^n\,w_iu_i$.  The {\em tropicalization\/}
$\tau(Y)$ of~$Y$ is the set
\begin{equation} \label{eq_trop}
\tau(Y)\,\, = \,\ \{w\in \R^n\,:\, {\rm in}_{w}(I_Y) 
                    \mbox{ does not contain a monomial}\,\}\, . 
\end{equation}

Let $K = \C\{ \!\{ \epsilon^\R \} \!\}$ be the field of Puiseux
series, i.e., series with complex coefficients, real exponents and well
ordered supports. The elements of
$K$ are also known as {\em transfinite Puiseux series}, and they form an
algebraically closed field of characteristic zero. The {\em order\/}
of a non-zero element $z$ in $K^* = K \backslash \{0\}$ is the
smallest real number $\nu$ such that $\epsilon^\nu$ appears with
non-zero coefficient in $z$. For a vector $z = (z_1,\ldots,z_s)$ in
$(K^*)^s$ we write $\,{\rm order}(z) := ({\rm order}(z_1), \ldots,
{\rm order}(z_s)) \in \R^s$. The points in the tropicalization
$\tau(Y)$ are precisely the orders of
  $K^*$-valued points on the variety $Y$ (see \cite{EKL},  \cite[Theorem 2.1.2]{Spe},
\cite[Theorem 2.1]{SS}). 

The set $\tau(Y)$ carries the structure of a polyhedral fan.
Namely, it is a subfan of the
Gr\"obner fan of $I_Y$; see \cite[\S 9]{S}. By a result of Bieri and Groves~\cite{BG}, 
the fan $\tau(Y)$ is pure of dimension $r$.
In \cite{BJSST} it was shown that $\tau(Y)$ is connected in codimension one,
and a practical algorithm was given for computing $\tau(Y)$
from polynomial generators of $I_Y$. We
will view the tropicalization  $\tau(Y)$ of a projective variety as an
$(r{-}1)$-dimensional fan in {\em tropical projective space\/} $\T
\P^{n-1}\,{:=}\,\R^n/\R(1,1,\ldots,1)$, which is an 
$(n{-}1)$-dimensional real affine space.

Every maximal cone $\sigma$ of the fan $\tau(Y)$ comes naturally with
an {\em intrinsic multiplicity} $m_\sigma$, which is a positive
integer. The integer $m_\sigma$ is computed as the sum of the
multiplicities of all monomial-free minimal associated primes of the
initial ideal ${\rm in}_w(I_Y)$ in $\C[x_1,\ldots,x_n]$, where $w$ is
in the relative interior of the cone $\sigma$.

\begin{remark}
\label{speyerdiss}
A geometric description
of the intrinsic multiplicity $m_\sigma$ arises from the beautiful interplay
of degenerations and compactifications discovered by Tevelev \cite{Tev} 
and studied by Speyer \cite[Chapter 2]{Spe} and Hacking (unpublished).
Let ${\bf X}$ denote the toric variety associated with the fan $\tau(Y)$.
Consider the intersection $Y_0$ of $Y$ with the 
dense torus $T$ in $\C\P^{n-1}$, and let
$\overline{Y_0}$ be the closure of $Y_0$ in ${\bf X}$.
By \cite[1.7, 2.5, and 2.7]{Tev},
the variety $\overline{Y_0}$ is complete and the 
multiplication map $\Psi: T \times \overline{Y_0}
\rightarrow {\bf X}$ is faithfully flat.
If follows that 
the intersection of $\overline{Y_0}$ with
a codimension $k$ orbit has codimension $k$ in
$\overline{Y_0}$.
In particular, the orbit ${\mathcal O}(\sigma)$
corresponding to a maximal cone $\sigma$ of
$\tau(Y)$ intersects $\overline{Y_0}$ in 
a zero-dimensional scheme $Z_\sigma$.
The intrinsic multiplicity $m_\sigma$
of the maximal cone $\sigma$ in the tropical variety $\tau(Y)$
is the length of $Z_\sigma$.
\end{remark}

We list three fundamental examples which will be important for our work.

\smallskip

\noindent
{\bf (1)} Let $Y$ be a {\bf hypersurface} in $\C\P^{n-1}$ defined by
an irreducible polynomial $f$ in $\C[x_1,\ldots,x_n]$. Then $\tau(Y)$
is the union of all codimension one cones in the normal fan of the
Newton polytope of $f$. The intrinsic multiplicity $m_\sigma$ of each
such cone~$\sigma$ is the lattice length of the corresponding edge of
the Newton polytope of $f$.

\noindent
{\bf (2)} Let $Y\,{=}\,X_A$ be the {\bf toric variety} defined by an
integer $d \times n$-matrix $A$ as above.  Its tropicalization
$\tau(X_A)$ is the linear space spanned by the rows of~$A$.

\noindent
{\bf (3)} Let $Y$ be a {\bf linear subspace} in $\C^n$
or in $\C \P^{n-1}$.
The tropicalization $\tau(Y)$ is the Bergman fan 
of the matroid associated with $Y$; see  \cite{AK, FS, S}
and (\ref{TLSrep1}) below.

In the last two families of examples, all the intrinsic multiplicities
$m_\sigma$ equal  $1$.
 
\smallskip

The tropicalization $\tau(Y)$ can be used to compute 
numerical invariants of~$Y$.
First, the dimension of $\tau(Y)$ coincides with the
dimension of $Y$.  In  Theorem \ref{thm:mult} below, we
express the multiplicities of the minimal primes in the initial monomial
ideals of $I_Y$  in terms of $\tau(Y)$. Equivalently,
we compute the algebraic cycle
of any initial monomial ideal ${\rm in}_w(I_Y)$.
This formula tells us the degree of the variety $Y$,
namely, the degree is the sum of
the multiplicities of the minimal
primes of $\, {\rm in}_w(I_Y)$.

Let $c\,{:=}\,n {-}r$ denote the codimension of
the irreducible projective variety~$Y$ in $\C\P^{n-1}$. Assume that $Y$
is not contained in a coordinate hyperplane, and let $I_Y$ be its
homogeneous prime ideal in $\C[x_1,\ldots,x_n]$. If~$w$ is a
generic vector in $\R^n$, the initial ideal ${\rm in}_w(I_Y)$ is a
monomial ideal of codimension~$c$.  Every minimal
prime over ${\rm in}_w(I_Y)$ is generated by a subset of $c$ of
the variables. We write $\, P_\tau \, = \,\langle \,x_i \,:\,i
\in \tau \,\rangle $ for the monomial prime ideal indexed by
the subset $\tau = \{\tau_1,\ldots,\tau_c\} \subset
\{1,2,\ldots,n\}$.

Assume that  the cone $ w+ \R_{> 0}\{e_{\tau_1},\ldots, e_{\tau_c}\}$
  meets the tropicalization $\tau(Y)$.   
We may suppose that the generic weight vector  $w \in \R^n$ 
satisfies that the image of $w$ in $\T \P^{n-1}$ does not lie in $\tau(Y)$ and that
the intersection of the cone $ w+ \R_{> 0}\{e_{\tau_1},\ldots, e_{\tau_c}\}$
with $\tau(Y)$ is finite and contained in the union of the relative 
interiors of its maximal cones. 
Let $\sigma$ be a maximal cone of the tropical variety and
\begin{equation} \label{eq:LL'}
\{v \} \, = \, (L+w) \cap L', 
\end{equation}
 where $ L= \R \{e_{\tau_1},\ldots, e_{\tau_c}\}$ and 
$L' =\R \sigma$ are the corresponding linear spaces, which are defined
over $\Q$.
We associate with  $\,v\,$ the 
{\it lattice multiplicity} of the intersection of $L$ and $L'$, which is
defined as the absolute value
of the determinant of any $n \times n$-matrix
whose columns consist of a 
$\Z$-basis of $\Z^n \cap L$ and a
$\Z$-basis of $\Z^n \cap L'$.

Here is the main result of this section.

\begin{theorem} \label{thm:mult} 
For $w \in \R^n$ a generic weight vector, a prime ideal
$P_\tau$ is associated to the initial monomial ideal ${\rm in}_w(I_Y)$ if and
  only if the cone $ w+ \R_{> 0}\{e_{\tau_1},\ldots, e_{\tau_c}\}$
  meets the tropicalization $\tau(Y)$.  The number of intersections,
each counted with its associated lattice multiplicity  times the intrinsic multiplicity,
is the multiplicity of the monomial ideal ${\rm in}_w(I_Y)$
  along the prime $P_\tau$.
   \end{theorem}

\begin{proof}
We work over the Puiseux series field $\,K = \C\{\!\{\epsilon^\R\}\!\}$,
we  write $K \P^{n-1}$ for the
$(n-1)$-dimensional projective space over the field $K$,
and we consider $Y$ as a subvariety of
 $K \P^{n-1}$.  
 The  translated variety $\,\epsilon^{-w} \cdot Y\,$ 
  is defined by the prime ideal
  $$
  \epsilon^w \cdot I_Y \quad = \quad \bigl\{\, f (\epsilon^{w_1}
  x_1,\ldots,\epsilon^{w_n} x_n ) \,: \,f \in I_Y \,\bigr\} \quad
  \subset \quad K[x_1,\ldots,x_n ].$$
Clearly, a point of the form 
 $w + u$ lies in $\tau(Y)$ if and only if there is a point
 in $\, \epsilon^{-w} \cdot \, Y\,$ with order $u$.

Let $L$ be a general linear
  subspace of dimension $c$ in $K\P^{n-1}$ which is defined over~$\C$. 
  The intersection $\, \epsilon^{-w} \cdot Y \,
  \cap \, L \,$ is a finite set of reduced points in $K \P^{n-1}$. 
  The number of these points is the degree $d$ of $Y$. 
There is a flat family over $K$ with generic fiber
$\, \epsilon^{-w} \cdot \, Y\,$ and special fiber the scheme $V({\rm in}_w(I_Y))$. 
Since $L$ is generic, the intersection of $L$ 
with this family is still
flat. Then, all points in the intersection $\, \epsilon^{-w} \cdot \, Y\,$
with $L$ are liftings of points in $\, V({\rm in}_w(I_Y)) \, \cap \, L\,$
(and all points can be lifted).   Moreover, 
the multiplicity of the $P_\tau$-primary component of the initial
ideal ${\rm in}_w(I_Y)$ equals the number of points 
in $\, \epsilon^{-w} \cdot \, Y\, \cap \, L \,$ of the form
  $$
  \theta \cdot \epsilon^u + \ldots \quad = \quad \bigl( \,\theta_1
  \epsilon^{u_1} + \ldots \,:\,\theta_2 \epsilon^{u_2} + \ldots\,:\,
  \cdots \,:\, \theta_n \epsilon^{u_n} + \ldots \,\bigr) ,$$
  where $\theta_k \in \C^*$ for all $k$, 
    $\,u_i = 0 \,$ for $i \not \in \tau$ and $u_j > 0 \,$ for $j \in
  \tau $, i.e., the number of points with values in the intersection
  of $\tau(Y)$ with
 the cone $\,w \,+ \R_{> 0}\{e_{\tau_1},\ldots, e_{\tau_c}\}$.

By our genericity assumption, each such intersection point $v$ lies on some
maximal cone $\sigma$ of $\tau(Y)$ and it is  counted
with its multiplicity, which is the product of the intrinsic multiplicity
$m_\sigma$ times the lattice multiplicity of the transversal intersection
of rational linear spaces  in (\ref{eq:LL'}). This product can be understood by
means of the flat family discussed in
Remark \ref{speyerdiss}. Namely, it follows from the
$T$-invariance of the multiplication map $\Psi$ that the
scheme-theoretic fiber of $\Psi$ over any point
of ${\mathcal O}(\sigma)$ is isomorphic to
$T' \times Z_\sigma$, where $T'$ is the stabilizer
of a point in ${\mathcal O}(\sigma)$. Since~${\bf X}$ is normal,
$T'$ is a torus $(\C^*)^{r-1}$. By \cite[1.7]{Tev},
$\,(\C^*)^{r-1} \times Z_\sigma$ is the intersection of~$T$
with the flat degeneration of $Y$ in $\C \P^{n-1}$
given by the one parameter subgroup of $T$ specified by 
the rational vector $w$. Our construction above 
amounts to computing the intersection of $\, (\C^*)^{r-1} \times Z_\sigma\,$
with the corresponding degeneration of $L$.
 The lattice index is obtained from the factor $(\C^*)^{r-1}$,
 and $m_\sigma$ is obtained from the factor $Z_\sigma$.
 Their product is the desired intersection number.
\end{proof}

The algebraic cycle of the variety $Y$ is represented by its {\em Chow
  form} $Ch_Y$, which is a polynomial in the bracket variables $[\gamma] =
[\gamma_1 \cdots \gamma_c]$; see \cite[\S 3.2.B]{GKZ}. 
Theorem~\ref{thm:mult} implies that the $w$-leading term of
the Chow form $Ch_Y$ equals
$\,\prod [\gamma]^{u_\gamma} $, where $u_\gamma$ is the (correctly
counted) number of points in $\,\tau(Y) \,\cap \, ( w+ \R_{>
  0}\{e_{\gamma_1},\ldots, e_{\gamma_c}\}) $.  We discuss this
statement for the three families of examples considered earlier.

\noindent
{\bf (1)} If $Y$ is a hypersurface then $c = 1$ and the
bracket variable $[\gamma] $ is simply the
ordinary variable $x_i$ for $i = \gamma_1$.
The $w$-leading monomial of the defining
irreducible polynomial $f(x_1,\ldots,x_n)$
equals $x_1^{u_1} \cdots x_n^{u_n}$
where $u_i$ is the number of times
the ray $w + \R_{>0} e_i$ intersects the
tropical hypersurface $\tau(Y)$, counted
with multiplicity.

\noindent
{\bf (2)} If $Y= X_A$ is a toric variety then $\tau(Y)
= {\rm rowspace}(A)$ and Theorem \ref{thm:mult} implies the familiar
result \cite[Thm.~8.3.3]{GKZ} that the initial cycles of $X_A$
are the regular triangulations of $A$. Indeed,  $ w+ \R_{>
  0}\{e_{\gamma_1},\ldots, e_{\gamma_{n-d}}\}$ intersects the row
space of $A$ if and only if the $(d-1)$-simplex $\bar \gamma = \{a_i :
i \not\in \gamma \}$ appears in the regular triangulation~$\Pi_w$ of
$A$ induced by~$w$. 
For a precise definition of $\Pi_w$ see Section~\ref{sec:subdivisions}. 
The intersection multiplicity is the lattice volume of $\bar \gamma$.
Hence $\,{\rm in}_w(Ch_{X_A}) = \prod_{\bar \gamma \in \Pi_w}
[\gamma]^{{\rm vol}(\bar \gamma )}$.

\noindent
{\bf (3)} If $Y$ is a linear space then
its ideal $I_Y$ is generated by $c$ linearly
independent linear forms in $\C[x_1,\ldots,x_n]$.
 The tropical variety $\tau(Y)$ is
the Bergman fan, to be discussed in Section~\ref{sec:monomials},
and the Gr\"obner fan of $I_Y$ is the normal fan of the associated
matroid polytope \cite{FS}.
For fixed generic $w$, there is a unique $c$-element subset
$\gamma$  of~$[n]$ such that
$ w+ \R_{> 0}\{e_{\gamma_1},\ldots, e_{\gamma_{c}}\}$ intersects $\tau(Y)$.
The intersection multiplicity is one,
and the corresponding initial ideal equals
$\,{\rm in}_w(I_Y) = \langle x_{\gamma_1} , \ldots, x_{\gamma_c} \rangle$.


\section{Tropicalizing maps defined by monomials in linear forms}
\label{sec:monomials}

In this section we examine a class of rational varieties $Y$
whose tropicalization $\tau(Y)$ can be computed combinatorially,
without knowing  the ideal $I_Y$.
We consider a rational map $f : \C^m \dashrightarrow \C^s$ that factors as
a linear map $\C^m \rightarrow \C^r$ followed by a Laurent monomial
map $\C^r \dashrightarrow \C^s$. The linear map is specified by a complex
$r\,{\times}\,m$-matrix $U = (u_{ij})$, and the Laurent monomial map
is specified by an integer $s\,{\times}\,r$-matrix $V = (v_{ij})$.
Thus the $i$-th coordinate of the rational map $f : \C^m \dashrightarrow
\C^s $ equals the following monomial in linear forms:
\begin{equation} \label{eq:linmonmap}
 f_i(x_1,\ldots,x_m)
\quad = \quad \prod_{k=1}^r \,(u_{k1} x_1 + \cdots + u_{km} x_m)^{v_{ik}}\,,
\quad  i=1,\ldots,s\, .
\end{equation}
Let $Y_{UV}$ denote the Zariski closure of the image of $f$. Observe
that if all row sums of $V$ are equal then
$f$ induces a rational map $\, \C \P^{m-1} \dashrightarrow \C \P^{s-1}$,
and the closure of its image is an irreducible projective variety, which we
also denote by $Y_{UV}$.

Our goal is to compute the tropicalization $\tau(Y_{UV})$ of the
variety $Y_{UV}$ in terms of the matrices $U$ and $V$.
In particular, we will avoid any reference to the equations of~$Y_{UV}$. 
The general framework of this section will be crucial for our proofs
of the results on $A$-discriminants and their tropicalization stated in 
the Introduction.

\smallskip
We list a number of special cases of varieties which have the form $Y_{UV}$.

\noindent
{\bf (1)}
If $r\,{=}\,s$, and $V\,{=}\,I_r$ then $f$ is the linear map $x \mapsto Ux$, 
and $Y_{UV} $ is the image of~$U$. We denote this
{\bf linear subspace} of $\C^r$ by ${\rm im}(U)$.
Its tropicalization $\tau({\rm im}(U))$ is the Bergman fan
of the matrix $U$, to be discussed in detail below.

\noindent
{\bf (2)} If $m\,{=}\,r$ and $U\,{=}\,I_m$ then $f$ is the monomial
map specified by the matrix $V$. Hence $Y_{UV}$ coincides with the
{\bf toric variety} $X_{V^t}$ which is associated with the transpose
$V^t$ of the matrix $V$. Its tropicalization is the column space
of~$V$.

\noindent
{\bf (3)}
Let $m\,{=}\,2$,
suppose the rows of $U$ are linearly independent,
and suppose the matrix $V$ has constant row sum.
Then $ Y_{UV} = {\rm image}( \C \P^1 \rightarrow \C \P^{s-1})$
is a {\bf rational curve}. Every rational projective curve 
arises from this construction, 
since every binary form is a monomial in linear forms.

Our next theorem implies that $\tau(Y_{UV})$ consists
of the rays in $\T \P^{s-1}$ spanned by the rows of $V$.
 Theorem \ref{thm:tUV} can also be derived from
 \cite[Proposition 3.1]{Tev}, but we prefer to give
 a proof that is entirely self-contained.

\begin{theorem} \label{thm:tUV}
The tropical variety $\tau(Y_{UV})$ equals the image
of the Bergman fan $\tau({\rm im}(U))$ under the linear map
$\R^r \rightarrow \R^s$ specified by the matrix $V$.
\end{theorem}

\begin{proof}
In what follows we consider all
algebraic varieties over the field $K = \C\{ \!\{ \epsilon^\R \}
\!\}$ and we use the characterization of the tropical variety $\tau(Y_{UV})$
in terms of $K^*$-valued points of the ideal of $Y_{UV}$. 
Extending scalars, we consider the map $f : K^m \dashrightarrow K^s$.
Let $z = f(x)$ be any point in the image of that map.  For $k \in
\{1,\ldots,r\}$ we set $ y_k = u_{k1} x_1 + \cdots + u_{km} x_m $
and $\alpha_k = {\rm order}(y_k)$. The vector $y = (y_1,\ldots,y_r)$
  lies in the linear space ${\rm im}(U)$, and hence the vector $\alpha =
  (\alpha_1,\ldots,\alpha_r)$ lies in    $\tau({\rm im}(U))$.
  
  The order of the vector $z = f(x) \in K^s $ is the vector $\, V
  \cdot \alpha \in \R^s$, since the order of its $i$th coordinate
  $z_i$ equals $\,\sum_{k=1}^r v_{jk} \cdot \alpha_k $.  Hence $z =
  f(x)$ lies in $V \cdot \tau({\rm im}(U))$, the image of the Bergman
  fan $\tau({\rm im}(U))$ under the linear map $V$.
  
  The image of $f$ is Zariski dense in $Y_{UV}$, i.e., there
  exists a proper subvariety  $Y$ of $Y_{UV}$ such that
  $Y_{UV} \backslash Y$ contains the image of $f$.
  By the Bieri-Groves Theorem~\cite{BG}, $\tau(Y)$
  is a polyhedral fan of lower dimension inside the
  pure-dimensional fan $\tau(Y_{UV})$. From this it follows that
  $\tau(Y_{UV})$ is the closure of $\tau(Y_{UV}\backslash Y)$ in
  $\R^s$. In the previous paragraph we showed that
  $\tau(Y_{UV}\backslash Y)$ is a subset of $V \cdot \tau({\rm im}(U))$.
  Since the latter is closed, we conclude that $\tau(Y_{UV}) \subseteq
  V \cdot \tau({\rm im}(U))$.
  
  It remains to show the converse inclusion $\,V \cdot \tau({\rm
    im}(U)) \subseteq \tau(Y_{UV})$.  Take any point $\beta \in V
  \cdot \tau({\rm im}(U))$, and choose $\alpha \in \tau({\rm im}(U))$
  such that $V \cdot \alpha = \beta$. There exists a $K^*$-valued
  point $y$ in the linear space ${\rm im}(U)$ such that ${\rm
    order}(y) = \alpha$. Then the point $\,y^V = \bigl(\prod_{k=1}^r
  y_k^{v_{1k}},\ldots, \prod_{k=1}^r y_k^{v_{sk}}\bigr) \,$ lies in
  $(K^*)^s \cap Y_{UV}$, and its order clearly equals $\beta$. Hence
  $\beta \in \tau(Y_{UV})$ as desired.
\end{proof}

\begin{remark} \label{intrmultYUV}
The intrinsic multiplicity $m_\sigma$ of any
maximal cone $\sigma$ in  a sufficiently fine fan structure on
the tropical variety $\tau(Y_{UV})$
is a lattice index which can be read off from the matrix~$V$.
Namely, suppose $\sigma$ 
lies the image of a maximal cone
$\sigma'$ of the Bergman fan $\tau({\rm im}(U))$. Then
$m_\sigma$ is the index of the subgroup 
$\,V( \R \sigma' \cap \Z^r)\,$ of 
$\R \sigma \cap \Z^s$.
This follows from Remark \ref{speyerdiss}
using standard arguments of toric geometry.
\end{remark}

Theorem \ref{thm:tUV} and Remark \ref{intrmultYUV}
furnish a combinatorial construction for
the tropicalization of any variety which is parameterized
by monomials in linear forms. Using
the results of Section~\ref{sec:tropical},
Theorem \ref{thm:tUV} can now be applied to
compute geometric invariants of such a variety,
for instance, its dimension, its degree and its initial cycles.
To make this computation effective,
we need an explicit description of the
Bergman fan $\tau({\rm im}(U))$.
Luckily, the relevant combinatorics is well understood, thanks to
\cite{AK, FS}, and in the remainder of this section
we summarize what is known.

Let $M$ denote the {\em matroid} associated
with the rows of the $r \times m$-matrix $U$.
Thus $M$ is a matroid of rank at most $m$
 on the ground set $[r] = \{1,2,\ldots,r\}$.
 The {\em bases} of~$M$ are the maximal subsets of $[r]$
 which index linearly independent rows of $U$.
Fix a vector $w \in \R^r$. Then the {\em $w$-weight} of
a basis $\beta $ of $M$ is $\sum_{i \in \beta} w_i$.
Consider the set of all bases of $M$
that have maximal $w$-weight.
This collection is the set of bases of a new matroid
which we denote by $M_w$. Note that each
$M_w$ has the same rank and the same ground set as $M = M_0$.
An element $i $ of $[r]$ is a {\em loop} of $M_w$ if
 it does not lie in any basis
of $M$ of maximal $w$-weight.

We can now describe 
$\tau({\rm im}(U))$ in terms of the matroid~$M$:
\begin{equation}
\label{TLSrep1}
\tau({\rm im}(U)) \quad = \quad
 \{w \in \R^r\,:\, M_{w} \mbox{ has no loop}\,\} \,.
\end{equation}
This representation endows our tropical linear space with the
structure of a polyhedral fan.  Namely, if $w \in \tau({\rm im}(U)) $,
then the set of all $w' \in \R^r$ such that $M_{w'} = M_{w}$ is a
relatively open convex polyhedral cone in $\R^r$.  The collection of
these cones is denoted $\BB(M)$ and is called the {\em Bergman fan} of
the matroid $M$. Depending on
the context, we may also write $\BB(U)$ and call it the Bergman fan of
the matrix $U$.


We now recall the connection between
 Bergman fans and {\em nested set complexes}.
 The latter encode the structure  of wonderful
compactifications of hyperplane arrangement complements in the work of
De~Concini and Procesi~\cite{DP}, and they were later studied from a
combinatorial point of view in~\cite{FK,FM,FY}.

A subset $X \subseteq [r]$ is a {\em flat} of our matroid $M$ if there
exists a vector $u \in {\rm im}(U)$ such that $X = \{i \in [r] : u_i
\,{=}\, 0 \}$. The set of all flats, ordered by inclusion, is the {\em
  geometric lattice} $\mathcal L \,{=}\,{\mathcal L}_M$.  A flat $X$
in ${\mathcal L}$ is called {\em irreducible\/} if the lower interval
$\{Y\,{\in}\,\mathcal L\,:\, Y\leq X\}$ does not decompose as a direct
product of posets.  Denote by $\mathcal I$ the set of irreducible
elements in $\mathcal L$. In other contexts, the irreducible elements
of a lattice of flats of a matroid were named {\em connected
  elements\/} or {\em dense edges\/}.  The matroid $M$ is {\em
  connected} if the top rank flat $\hat 1\,{=}\,[r]$ is irreducible,
and we assume that this is the case. Otherwise, we artificially add
$\hat 1$ to $\mathcal I$. We call a subset $S \subseteq \mathcal I$
{\em nested\/} if for any set of pairwise incomparable elements
$X_1,\ldots, X_t$ in $S$, with $t \geq 2$, the join $X_1\vee \ldots
\vee X_t$ is not contained in~$\mathcal I$. The nested subsets in
$\mathcal I$ form a simplicial complex, the {\em nested set complex\/}
${\mathcal N}(\mathcal L)$. See~\cite[Sect.\ 2.3]{FK}
for further information.

Feichtner and Yuzvinsky \cite[Eqn (13)]{FY} introduced the following
natural geometric realization of the nested set complex
${\mathcal N}(\mathcal L)$.  Namely, the collection of cones
\begin{equation}\label{eq:nscones}
       \R_{\geq 0} \{\,e_X\,:\, X\in S\,\} \qquad \mbox{ for }\,
       S\in {\mathcal N}(\mathcal L) 
\end{equation}
forms a unimodular fan whose
face poset is the face poset of the nested set complex of $\mathcal L$. 
Here $\,e_X = \sum_{i \in X} e_i\,$ denotes the
incidence vector of a flat $X \in \mathcal I$.
We consider this fan in the 
tropical projective space $\T \P^{n-1}$,
and we also denote it by ${\mathcal N}(\mathcal L)$.

It was shown in \cite[Thm 4.1]{FS} that
the {\em nested set fan} $\mathcal N(\mathcal L_M)$ is a simplicial
subdivision of the Bergman fan ${\mathcal B}(M)$, and hence
of our tropical linear space $\tau({\rm im}(U))$.
The Bergman fan need not be simplicial, so 
the nested set fan can be finer than the Bergman fan.
However, in many important cases the two fans
coincide~\cite[\S 5]{FS}.

What we defined above is the coarsest in a hierarchy of nested set
complexes associated with the geometric lattice $\mathcal L $. Namely,
for certain choices of subsets ${\mathcal G}$ in~$\mathcal L$, the
same construction gives a nested set complex ${\mathcal N}({\mathcal
L},{\mathcal G})$ which is also realized as a unimodular simplicial
fan.  Such ${\mathcal G}$ are called {\em building sets\/}; there is
one nested set fan for each building set ${\mathcal G}$ in $\LL$,
see~\cite[Sect.\ 2.2, 2.3]{FK}.  If two building sets are contained in
another, ${\mathcal G}_1\,{\subseteq}\ {\mathcal G}_2$, then
${\mathcal N}({\mathcal L},{\mathcal G}_2)$ is obtained from
${\mathcal N}({\mathcal L},{\mathcal G}_1)$ by a sequence of stellar
subdivisions~\cite[Thm 4.2]{FM}.  The smallest building set is
${\mathcal G} = {\mathcal I}$, the case discussed above, and the
largest building set is the set of all flats, ${\mathcal G} =
{\mathcal L}$.  In the letter case, the corresponding nested set
complex ${\mathcal N}({\mathcal L},{\mathcal L})$ is the {\em order
complex\/} or {\em flag complex\/} 
$\mathcal F({\mathcal L})$ of the lattice, i.e., the simplicial
complex on ${\mathcal L}\,{\setminus}\,\{\hat 0\}$ whose simplices are
the totally ordered subsets.  Again, there is a realization of
$\mathcal F({\mathcal L})$ as a unimodular fan with rays generated by the
incidence vectors of flats in~$\mathcal L$. We call this fan the {\em
flag fan\/} of~$\mathcal L$ and denote it by $\mathcal F({\mathcal L})$
for ease of notation.

Summarizing the situation, we have the following sequence
of subdivisions each of which can be used to compute
the tropicalization of a linear space. 

\begin{theorem}
\label{threeSubdivisions}
Given a matrix $U$ and the matroid $M$ of rows in~$U$,
the tropical linear space $\tau({\rm im}(U))$ has
 three natural fan structures: the Bergman fan
${\mathcal B}(M)$ is refined by the
nested set fan ${\mathcal N}({\mathcal L}_M)$, which  is
refined by the flag fan $\mathcal F({\mathcal L}_M)$.
\end{theorem}

We present an example which illustrates the concepts 
developed in this section.

\begin{example} \label{ex:UV}
Let $m{=}4$, $r{=}5$, $s{=}4$, and consider the map~$f{=}(f_1,f_2,f_3,f_4)$
from $\C^4$ to $\C^4$ 
whose coordinates are the
following monomials in linear forms:
\begin{eqnarray*} 
& f_1 \quad =&  (x_1-x_2)^3 (x_1-x_3)^3 \\
& f_2 \quad =& (x_1-x_2)^2 (x_2-x_3)^2 (x_2-x_4)^2 \\
&  f_3 \quad =& (x_1-x_3)^2 (x_2-x_3)^2 (x_3-x_4)^2 \\
&  f_4 \quad = & (x_2-x_4)^3 (x_3-x_4)^3\, .
\end{eqnarray*}
So $f$ is as in (\ref{eq:linmonmap}) with
$$ U \,\,\, = \,\,\,
\left(
\begin{array} {rrrr}
               1  & -1 &  0 &  0  \\
                1  &  0 & -1 &  0  \\
                0  &  1 & -1 &  0  \\ 
                0  &  1 &  0 & -1  \\ 
                0  &  0 &  1 & -1  
 \end{array}  \right) \quad \hbox{and} \quad
 V \,\,\,= \,\,\,
 \begin{pmatrix}
 3 & 3 & 0 & 0 & 0\\
 2 & 0 & 2 & 2 & 0\\
 0 & 2 & 2 & 0 & 2\\
 0 & 0 & 0 & 3 & 3
 \end{pmatrix}. $$ 
 The projectivization of the variety $Y_{UV}{=}{\rm cl}({\rm image}(f))$
 is a surface in $\C \P^3$, since $U$ has rank $3$. We construct
the irreducible homogeneous polynomial $P(z_1,z_2,z_3,z_4)$ 
which defines this surface.
The matroid of~$U$ is the graphic matroid of $K_4$ with one edge removed.
{}From \cite[Example 3.4]{FS} we know that
the Bergman complex $\BB(U)$ is 
 the complete bipartite graph~$K_{3,3}$, embedded
 as a $2$-dimensional fan in $\T\P^{4}$.
 
 The tropical surface $\tau(Y_{UV})$ is the image of $\BB(U)$ under
 $V$.  This image is a $2$-dimensional fan in $\T\P^{3}$.
 It has seven rays: six of them are images of the
 rays of $\BB(U)$, the last one  is the intersection of the images of two
 $2$-dimensional faces that occurs due to the non-planarity of~$K_{3,3}$. 
 Hence the Newton polytope of the polynomial $P$ is $3$-dimensional
 with $6$ vertices, $11$ edges, and $7$ facets; see Figure~\ref{fig:NP}.  
 
 The six extreme monomials of $P$ can be computed (even by hand) using
 Theorem~\ref{thm:mult}, namely, by intersecting the rays
 $w+\R_{>0}e_i$ with $\tau(Y_{UV})$ in $\T\P^3$.  This computation
 reveals in particular that the degree of the polynomial $P$ is~$28$.
 Using linear algebra, it is now easy to determine all $171$ monomials
 in the expansion of~$P$.

\begin{figure}[ht]
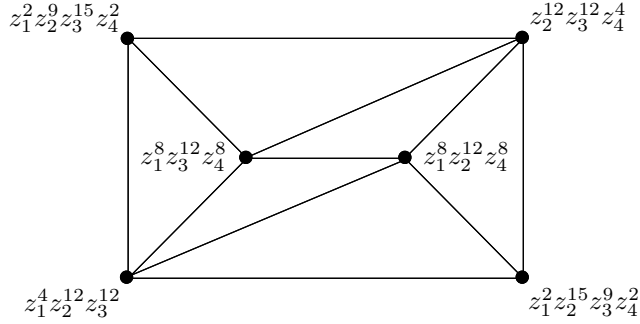

   \begin{picture}(0,0)%
     \includegraphics{UV.pstex}%
   \end{picture}%
   \input{UV.pstex_t}%
   
\caption{The Newton polytope of the polynomial $P$ in Example~\ref{ex:UV}.} 
\label{fig:NP}
\end{figure}
\end{example}


\section{Back to $A$-discriminants}
\label{sec:Horn}

In this section we return to the setting of the Introduction, and we
prove Theorem~\ref{thm:tropical} and Theorem
\ref{thm:initial}. Recall that $A$ is an integer
$d\,{\times}\,n$-matrix such that $(1,\ldots,1)$ is in the row span
of~$A$, i.e., the column vectors $a_1,a_2,\ldots,a_n$ lie in an affine
hyperplane in~$\R^d$. We also assumed that the vectors $a_1,\ldots,a_n$
span~$\Z^d$. We identify the matrix $A$ with the point configuration
$\{a_1,a_2,\ldots,a_n\}$.  The convex hull of the configuration $A$ is
a $(d{-}1)$-dimensional polytope with $\leq n$ vertices.

The projective toric variety $X_A$ is defined as the closure of the
image of the monomial map $\psi_A : (\C^*)^d \to \C \P^{n-1},\,$ $\,t
\mapsto (t^{a_1}: t^{a_2} \, : \, \cdots \, : \, t^{a_n})$.
Equivalently, $X_A$ is the set of all points $x \in \C \P^{n-1}$ such
that $x^u = x^v$ for all $u,v \in \N^n$ with $Au = Av$.

Let $(\C \P^{n-1})^* $ denote the projective space dual to $\C
\P^{n-1}$. The point $\xi = (\xi_1 : \dots : \xi_n)$ in $\, (\C
\P^{n-1})^* \,$ corresponds to the hyperplane $H_\xi =\{ x \in \C
\P^{n-1} \, : \, \sum_{i=1}^n x_i \xi_i =0\}$. The {\em dual variety}
$X_A^*$ is defined as the closure in $(\C \P^{n-1})^*$ of the set of
points~$\xi$ such that the hyperplane $H_\xi$ intersects the toric
variety $X_A$ at a regular point $p$ and contains the tangent space
$T_{X_A}(p)$ of $X_A$ at $p$.

Kapranov \cite{K} showed that reduced discriminantal varieties are parametrized
by monomials in linear forms. This parametrization, called the
{\em Horn uniformization}, will allow us to determine the tropical
discriminant $\tau (X_A^*)$ via the results of Section~\ref{sec:monomials}.
We denote by $\C\P({\rm ker}(A))$ the projectivization of the kernel
of the linear map given by~$A$, an $(n{-}d{-}1)$-dimensional
projective subspace of $\C \P^{n-1}$, and we denote by $T^{d-1}
= (\C^*)^d/\C^*$ the dense torus of $X_A$.
The following result is a variant of 
\cite[Theorem 2.1]{K}; see also \cite[\S 9.3.C]{GKZ}.

\begin{proposition} \label{prop:HDparametrization}
 The dual variety $X_A^*$ of the toric variety $X_A$ is the closure of
  the image of the map $ \, \varphi_A : \C \P({\rm ker}(A))\times
  T^{d-1} \to (\C \P^{n-1} )^* \,$ which is given by
\begin{equation} \label{eq:phiA}
 \varphi_A(u,t) \,\,=\,\, (u_1 t^{a_1} \, :  \,u_2 t^{a_2}
 \, :\, \cdots \, : \, u_n t^{a_n}). 
\end{equation}
\end{proposition}

\begin{proof}
  Consider the unit point ${\bf 1} = (1:1:\ldots:1)$ on the toric
  variety~$X_A$.  The hyperplane $H_\xi$ contains both the point ${\bf
    1}$ and the tangent space $T_{X_A}({\bf 1})$ at this point if and
  only if $\xi $ lies in the kernel of $A$. This follows by evaluating
  the derivative of the parametrization $\psi_A$ of $X_A$ at
  $\,(t_1,t_2,\ldots,t_d) = (1,1,\ldots,1)$.  If $p = \psi_A(t)$ is
  any point in the dense torus of $X_A$, then the tangent space at
  that point is gotten by translating the tangent space at ${\bf 1}$
  as follows:
  $$T_{X_A}(p) \quad = \quad p \cdot T_{X_A}({\bf 1}) .$$
  The
  hyperplane $H_\xi$ contains $p$ if and only if $p^{-1} \cdot H_\xi =
  H_{\xi \cdot p}$ contains ${\bf 1}$, and~$H_\xi$ contains
  $T_{X_A}(p)$ if and only if $ H_{\xi \cdot p}$ contains
  $T_{X_A}({\bf 1})$.  These two conditions hold, for some $p $ in the
  dense torus of $X_A$, if and only if $\,\xi \in {\rm
    image}(\varphi_A)$.
\end{proof}

Proposition \ref{prop:HDparametrization} shows that the dual variety
$X_A^*$ of the toric variety $X_A$ is parametrized by monomials in
linear forms. In the notation of Section~\ref{sec:monomials} we set $m
= n$, $r = n+d$, $s = n$, and the two matrices are
\begin{equation}
\label{chooseUandV}
U \,\,=\,\, \begin{pmatrix} \, B & 0  \\ 0 & I_d \, \end{pmatrix} 
\quad \hbox{and} 
\quad
V \,\,= \,\, \begin{pmatrix} \, I_n & A^t  \, \end{pmatrix}, 
\end{equation}
where $B$ is an $n\times (n{-}d)$-matrix whose columns
span  the kernel of $A$ over the integers.
Thus the rows of $B$ are {\em Gale dual} to the configuration $A$.

\begin{lemma} \label{lem:connectToSec3}
The variety $Y_{UV} \subset \C^n $ defined by {\rm (\ref{chooseUandV})} 
as in Section~\ref{sec:monomials}
is equal to the cone over  the dual variety $X_A^*  \subset
(\C \P^{n-1})^* $ of the toric variety $X_A
 \subset \C \P^{n-1}$.
\end{lemma}

\begin{proof}
  Let $f $ be the rational map defined by (\ref{chooseUandV}) as in
  Section~\ref{sec:monomials}, and set $x = (x_1,\ldots,x_{n-d})$ and
  $t = (t_1,\ldots,t_d) = (x_{n-d+1}, \ldots, x_n) $.  Then
  (\ref{eq:linmonmap}) equals
$$ f_i(x,t) \quad = \quad
(b_{i,1} x_1 + \cdots + b_{i,n-d} x_{n-d}) \cdot t^{a_i}, $$
which equals the $i$-th coordinate of $\varphi_A$
if we write ${\rm ker}(A)$ as the image of $B$.
\end{proof}

We are now ready to prove Theorem \ref{thm:tropical}.

\smallskip

\noindent
{\em Proof of Theorem~\ref{thm:tropical}.\/}
Let us first note that the {\em co-Bergman fan} 
${\mathcal B}^*(A)$ of the rank~$d$ configuration given by the columns 
of $A$ equals the Bergman fan of the rank  $n-d$ configuration given by 
the rows of $B$.
Thus ${\mathcal B}^*(A)$, as it appears in~(\ref{eq:coBergmanFan}) 
is the Bergman fan of the matroid dual to the matroid given by the columns 
of~$A$.

The support of the co-Bergman fan ${\mathcal B}^*(A)$ 
is the tropicalization of the linear space
$\,{\rm ker}(A) = {\rm im}(B)$.
Now, if $U$ is taken as in (\ref{chooseUandV}) then 
we have the following decomposition in $\,\R^r = \R^n \oplus \R^d$:
$$ \tau({\rm im}(U)) \quad = \quad
\tau({\rm im}(B)) \oplus \tau({\rm im}(I_d)) \quad = \quad
{\mathcal B}^*(A) \oplus \R^d. $$
The image of this fan under the linear map
$\,V \,= \, \begin{pmatrix} \, I_n \! & \! A^t \, \end{pmatrix}\,$
is the (Minkowski) sum of  ${\mathcal B}^*(A)$ and the
image of $A^t$. Of course, the latter is the
row space of $A$. Hence our assertion
follows from Lemma \ref{lem:connectToSec3} and
Theorem~\ref{thm:tUV}. \hfill $\Box$

\medskip
Similarly, the tropicalization of the reduced version
of the dual variety can be derived from Theorem~\ref{thm:tUV}. Let
 $B$  be again an $n\times(n-d)$-matrix whose columns span the
kernel of~$A$ over the integers. The
{\em reduced dual variety\/} $Y_B^*:= Y_{BB^t}\,$ is the closure of the image of the
rational morphism $\widetilde \varphi_B : \C\P^{n-d-1} \dashrightarrow
  \C^{n-d}$, whose
$i$-th coordinate equals
\[
\widetilde \varphi_B(s_1:\ldots:s_{n-d})_i\, \ = \, \,
\prod_{k=1}^n\,(\,b_{k,1}s_1 + \ldots + b_{k,n-d}s_{n-d}\,)^{b_{k,i}}\, ,
\quad i=1,\ldots, n-d\, .
\]

\begin{corollary} \label{crl:reducedDualVariety}
Given a Gale dual $B$ of~$A$, the
tropicalization of the reduced dual variety $Y_B^*$ is the
image of the co-Bergman fan ${\mathcal B}^*(A)$ under the linear map $B^t$. 
\end{corollary}

We often do not distinguish between the reduced and unreduced version
of dual varieties and their tropicalizations, and denote both by $X_A^*$
and $\tau(X_A^*)$, respectively.

\smallskip

We illustrate Theorem~\ref{thm:tropical} for the case
when $X_A $ is the Veronese surface, regarded as the
projectivization of the variety
of all symmetric $3 \times 3$ matrices of rank $\leq 1$.

\begin{example} \label{ex:symm3x3}
We take $d=3, n = 6$, and we fix the matrix
$$ A \quad = \quad 
\begin{pmatrix}
1 & 1 & 1 & 1 & 1 & 1\\
0 & 1 & 2 & 0 & 1 & 0 \\
0 & 0 & 0 & 1 & 1 & 2 \\
\end{pmatrix}. $$
Note that the more usual matrix
$$ A' \,= \, 
\begin{pmatrix}
2 & 1 & 0 & 1 & 0 & 0\\
0 & 1 & 2 & 0 & 1 & 0 \\
0 & 0 & 0 & 1 & 1 & 2 \\
\end{pmatrix} \, = \, 
\left( \begin{array}{rrr}
2  & -1 & -1 \\
0 & 1 & 0 \\
0 & 0 & 1 
\end{array} \right) \, \cdot \, A,$$
defines the same Veronese embedding of
$\C\P^2$ into $\C\P^5$, but 
the columns of $A'$ do not span $\Z^3$ and the 
monomial para\-metri\-za\-tion $\psi_{A'}$
is two-to-one.

Points $x$ in $\C\P^5$ are identified with symmetric $3 \times 3$-matrices
$$ X \quad = \quad
\begin{pmatrix}
2 x_1 & x_2 & x_4 \\
   x_2 & 2 x_3 & x_5 \\
   x_4 & x_5   & 2 x_6 \\
\end{pmatrix}\, .
$$
A point $u$ is in $\C\P({\rm ker}(A))$ if and only if the
corresponding matrix $U$ has zero row  and column sums.
If this holds, and $t $ is any point in $ (\C^*)^3$, then
the symmetric $3 \times 3$-matrix $X$ corresponding to $x = \varphi_A(u,t)$ 
is singular because it satisfies
$$  \begin{pmatrix} 1 & 1/t_2 & 1/t_3 \end{pmatrix}  \cdot X 
\quad = \quad \begin{pmatrix} 0 & 0 & 0 \end{pmatrix}. $$
Hence $\varphi_A$  parametrizes  rationally the hypersurface of singular
symmetric matrices~$X$, and the
 $A$-discriminant equals the classical discriminant
 $\D_A (x) = \frac{1}{2}{\rm det}(X)$.

\begin{figure}[ht]
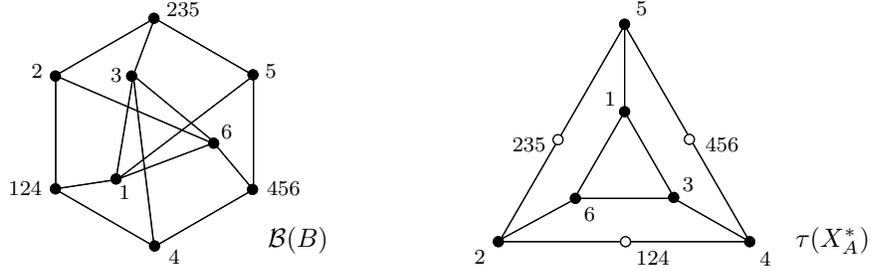

   \begin{picture}(0,0)%
     \includegraphics{B3x3.pstex}%
   \end{picture}%
   \input{B3x3.pstex_t}%
   
\caption{Bergman complex and tropical discriminant in Example~\ref{ex:symm3x3}}\label{fig_B3x3}
\end{figure}

The tropicalization of $X_A^*$ is obtained as follows.
We choose a Gale dual $B$ of~$A$,
$$   B^t \quad = \quad
\left( \begin{array}{rrrrrr}
  1 & -2 & 1 & 0 & 0 & 0 \\
  1 & -1 & 0  & -1 & 1 & 0 \\
  1 & 0& 0  & -2 & 0 & 1 
\end{array}
\right)
 \, .$$
Note that the matroid given by the columns of $A$ is self-dual. 
The Bergman fan $\BB(B)\,{=}\,\BB^*(A)$  
is a $2$-dimensional fan in $\R^6/ \R(1,1,1,1,1,1)$,
or, equivalently, a graph on the $4$-sphere. We depict this graph in
Figure~\ref{fig_B3x3} on the left. It has
nine vertices, corresponding to the six singletons $1,2,\ldots,6$ and the three
circuits $124$, $235$, $456$.

Its image under $B^t$ is the tropical discriminant $\tau(X_A^*)$, a
$2$-dimensional fan in $3$-dimensional real space or, equivalently, a
graph on the $2$-sphere.  We depict this graph in
Figure~\ref{fig_B3x3} on the right.  
The rays $124$, $235$, and $456$ of $\BB(B)$
map into the relative interiors of the $2$-dimensional
cones $24$, $25$, and $45$ of $\tau(X_A^*)$, respectively. More precisely,
the ray $124$  of $\BB(B)$ is mapped to the negative of the image of 
the ray $1$, which in turn is the dividing ray in the $2$-dimensional 
cone spanned by the images of the rays $2$ and $4$. 
The rays $235$ and $456$ are mapped similarly.
  
Since $A$ is non-defective, the tropical discriminant $\tau(X_A^*)$ is
the union of codimension~$1$ cones in the normal fan of the Newton
polytope of~$\Delta_A$. {}As $\tau(X_A^*)$ is the $1$-skeleton of
a triangular prism, we conclude that the Newton polytope of $\Delta_A$
is a bipyramid. Its five vertices correspond to the five terms in
the determinant of $X$.
\end{example}

\smallskip

Returning to the general case,
we note that the dimension of the image of $\varphi_A$ is at most $\, {\rm
  dim}(\C\P({\rm ker}(A)) \times T^{d-1}) \, = n-2 $, so the dual 
variety $X_A^*$ is a
proper subvariety of $\C\P^{n-1}$.  If the dimension of $X_A^*$ is less
than $n-2$, that is, if $X_A^*$ is not a hypersurface, we say that
the toric variety $X_A$ and its point configuration~$A$, 
are {\em defective}. In the non-defective case, 
there is a unique (up to sign) irreducible polynomial
$\D_A$ with integer coefficients
 which vanishes on $X_A^*$. The polynomial $\D_A$ is 
the {\em $A$-discriminant\/} as defined
 in \cite[\S 9.1.A]{GKZ}. In what follows, 
the dual variety $X_A^*$ itself will be referred
to as the $A$-discriminant, even if $A$ is defective.

By the Bieri-Groves Theorem \cite{BG, S}, the dimension
of the $A$-discriminant $X_A^*$ coincides with the dimension
of the tropical $A$-discriminant $\tau(X_A^*)$. 
Theorem \ref{thm:tropical} 
furnishes a purely combinatorial formula for that dimension.

\begin{corollary} \label{cor:dimension}
The dimension of the $A$-discriminant $X_A^*$ in 
$\C\P^{n-1}$ is one less than
the largest rank of any matrix
$ (A^t \!,\sigma_1, \dots, \sigma_{n-d-1})$
where $\sigma $ runs over $\,{\mathcal C}(A)$.
\end{corollary}

Here ${\mathcal C}(A)$ is the subset of $\{0,1\}^n$ defined in the
Introduction. That definition is now best understood using the
matroid-theoretic concepts which we reviewed in the second half of
Section~\ref{sec:monomials}, where we take $U $ to be the $n \times
(n-d)$-matrix $B$ as in~(\ref{chooseUandV}) and, hence, $M$ to be the
rank~$n{-}d$ matroid associated with the rows of~$B$.  In fact, $M$ is
the matroid dual to the matroid given by the columns of the $d \times
n$-matrix $A$, and $\LL(A)$ coincides with the lattice of
flats~$\LL_M$. The set ${\mathcal C}(A)$ corresponds to the facets of
the flag complex $\mathcal F({\mathcal L}_M)$. In light of Theorem
\ref{threeSubdivisions}, one could reformulate Corollary
\ref{cor:dimension} with $\sigma$ ranging over the facets of the
nested set complex $\mathcal{N}(\mathcal{L}_M)$ or the Bergman fan $
{\mathcal B}(M) = {\mathcal B}^*(A)$.

We are now prepared to state
and prove the general version of Theorem~\ref{thm:initial}.
Let ${\mathcal C}^c$ denote the set of all
proper chains of length $n-d-c-1$ in
${\mathcal L}(A) = {\mathcal L}_M$, where we identify flats of the matroid
in ${\mathcal L}(A)$ with their incidence 
vectors, and hence the chains in ${\mathcal C}^c$ with 
$(n{-}d{-}c)$-element subsets of $\{0,1\}^n$. Equivalently,
${\mathcal C}^c$ is the set of
$(n{-}d{-}c)$-element subsets of the elements of
${\mathcal C} = {\mathcal C}(A)$. 
We write ${\rm in}_w(X_A^*) = 
{\rm in}_w(I_{X_A^*})$ for the initial
ideal, with respect to some $w \in \R^n$, of the
homogeneous prime ideal $I_{X_A^*}$
of the $A$-discriminant $X_A^*$. 

\begin{theorem} \label{thm:initial2}
  Suppose that the $A$-discriminant $X_A^*$ has codimension $c$ and
  let $\tau = \{\tau_1 ,\ldots,\tau_c\} \subset \{1,\ldots,n\}$.  If
  $w$ is a generic vector in $\R^n$ then the multiplicity of the
  initial monomial ideal $\,{\rm in}_w(X_A^*)\,$ along the prime 
$\, P_\tau \,= \,\langle \,x_i \,:\,i \in \tau \,\rangle \,$  equals
\begin{equation}\label{eq:degree2}
\sum_{\sigma \in {\mathcal C}^c_{i,w} } |\, \det ( A^t,\sigma_1, \dots,
\sigma_{n-d-c}, e_{\tau_1}, \dots, e_{\tau_c} ) \,| \,  ,
\end{equation}
where ${\mathcal C}^c_{i,w}$ is the subset of ${\mathcal C}^c$
consisting of all chains $\sigma$ such that
\begin{equation}
\label{theyIntersect}
 {\rm rowspace}(A) \,\cap \,
 \R_{> 0} \bigl\{\sigma_1, \ldots, \sigma_{n-d-c}, 
-e_{\tau_1},\ldots, -e_{\tau_c}, -w\bigr\} \,\not= \,\emptyset .
\end{equation}
\end{theorem}

Theorem \ref{thm:initial} is the special case of Theorem
\ref{thm:initial2} when $A$ is non-defective, i.e., $c\,{=}\,1$ 
and $I_{X_A^*} $ is the principal
ideal generated by $ \Delta_A $. In that case, the initial monomial
ideal ${\rm in}_w(X_A^*)$ is generated by the
initial monomial ${\rm in}_w(\Delta_A)$.

\smallskip
\noindent
{\em Proof of Theorem \ref{thm:initial2}.\/} According
to Theorem~\ref{thm:mult}, the prime $P_\tau$ is associated to ${\rm
  in}_w(X_A^*)$ if and only if the polyhedral 
cone $ w+ \R_{>0}\{e_{\tau_1},\ldots, e_{\tau_c}\}$ meets the
tropicalization~$\tau(X_A^*)$, which was described 
in~Theorem~\ref{thm:tropical} as 
$\,{\mathcal B}^*(A) + {\rm rowspace}(A)$. 

The collection of cones $\R_{\geq 0} \sigma$ for $\sigma \in {\mathcal
C}$ forms a unimodular triangulation of the co-Bergman fan ${\mathcal
B}^*(A)$. This was proved by Ardila and Klivans~\cite{AK}, and we
discussed it in Theorem~\ref{threeSubdivisions}, calling $(\R_{\geq 0}
\sigma)_{\sigma \in {\mathcal C}}$ the flag fan $\mathcal F({\mathcal L}_M)$ 
of the matroid given
by the rows of a Gale dual $B$ of $A$. Therefore,
(\ref{theyIntersect}) characterizes when $ w+ \R_{>
0}\{e_{\tau_1},\ldots, e_{\tau_c}\}$ meets $\,\R_{\geq 0} \sigma +
{\rm rowspace}(A)\,$ for some $\sigma\in{\mathcal C}^c$. The
multiplicity of this intersection is precisely the stated $n \times
n$-determinant.  This can be derived from Remarks \ref{speyerdiss} and
\ref{intrmultYUV}.
\hfill \qed 

\medskip

Our degree formula for the $A$-discriminant can now be rephrased in
the following manner which is more conceptual and geometric.

\begin{corollary} \label{cor:mult}
A monomial prime $P_\tau$ is associated to  ${\rm
    in}_w(X_A^*)$ if and only if
    the cone $ w+ \R_{> 0}\{e_{\tau_1},\ldots, e_{\tau_c}\}$
meets the fan $\,{\mathcal B}^*(A) + {\rm rowspace}(A)$.
 The number of intersections, counted with
  multiplicity, is the multiplicity of ${\rm in}_w(X_A^*)$ along~$P_\tau$.
   \end{corollary}


\section{Computations, subdivisions, and singular 
tropical hypersurfaces}
\label{sec:subdivisions}

We start this section with a brief discussion of computational issues.
The formula for the extremal terms of $\Delta_A$ in
Theorem~\ref{thm:initial} gives rise to a practical method for computing
the Newton polytope of the $A$-discriminant $\Delta_A$. Indeed,
the co-Bergman fan of the matrix $A$ can be computed efficiently by gluing
local Bergman fans, as explained in \cite[Algorithm 5.5]{FS}.
See Examples 5.7, 5.8 and 5.9 in \cite{FS} for some non-trivial
computations. Extending the software used for those
computations, we wrote a {\tt maple} program for evaluating
the formula (\ref{eq:degree1}) in Theorem \ref{thm:initial}.
The input for our program consists of three positive integers $d,n,R$,
and a $d \times n$-matrix $A$ which is assumed to be non-defective.
The output is a list of initial monomials ${\rm in}_w(\Delta_A)$
of the $A$-discriminant $\Delta_A$, for $R$
randomly chosen vectors $w$ in $\N^n$.  
Our {\tt maple} implementation is available upon request from any of the  
authors.

Note that in case $A$ is non-defective, it is possible to recover 
$\Delta _A$ up to constant from its Newton polytope $N(\Delta_A)$
by solving a linear system of equations. Namely, consider a generic polynomial
$g$ with exponents in the lattice points of $N(\Delta_A)$. Imposing
the condition that $g$ vanishes on the image of the rational
parametrization of $\Delta_A$ given in Proposition~\ref{prop:HDparametrization} 
and Lemma~\ref{lem:connectToSec3} translates into a system of linear equations in the 
coefficients of $g$ whose solution space is one-dimensional.

\begin{example}
\label{mixeddiscriminant}
Let $d=4$ and $n=8$, and consider the matrix $A$ below. The corresponding
$A$-discriminant is the \emph{mixed discriminant} of the two bivariate
polynomials $f_1, f_2$ whose exponent vectors are read respectively
from the first four and last
four columns of the last two rows of the matrix. The non-vanishing of
$\Delta_A$ expresses the condition that the intersection $f_1=f_2=0$
is transversal.
What follows is the output of our {\tt maple} program on this input.
On a fast workstation, our code takes about half a second to compute 
the co-Bergman fan. Afterwards
it takes about one second per initial monomial.
So the total running time for this matrix is about $R$ seconds, where
$R$ is the number of iterations specified by the user:

\bigskip 
\begin{smaller}
\begin{verbatim}
                   [ 1     1     1     1     0     0     0     0]
                   [                                            ]
                   [ 0     0     0     0     1     1     1     1]
              A := [                                            ]
                   [ 2     3     5     7    11    13    17    19]
                   [                                            ]
                   [19    17    13    11     7     5     3     2]

                    Computing the Co-Bergman fan of A....
                         DONE. Time elapsed =  0.340

           The number of maximal cones in the Co-Bergman fan is  57
     48 of these cones map to codimension one in the tropical discriminant.
     What follows are  3  pairs of weight vectors and initial monomials:

                                                    28   35   35   28
         [446, 773, 680, 37, 925, 963, 765, 380], x1   x4   x5   x8

                                                   34   29   2   39   22
       [439, 464, 454, 360, 303, 279, 591, 583], x1   x4   x5  x6   x8

                                                 2   39   22   22   39   2
     [801, 447, 685, 447, 765, 775, 358, 498], x1  x2   x4   x5   x7   x8
     
\end{verbatim}

\end{smaller}
{}From this output we see that this $A$-discriminant is a polynomial of degree
$\,126$. We discuss the geometric meaning of this example 
in Section~\ref{sec:resultants}.
\end{example}

\medskip

Gel'fand, Kapranov and Zelevinsky \cite{GKZ} established the relationship
between the $A$-discriminant and the {\em secondary fan} of a point 
configuration. The secondary fan 
parametrizes the regular polyhedral subdivisions of $A$. It is shown in 
\cite{GKZ} that, in the non-defective case, the A-discriminant divides
the principal $A$-determinant. Hence, the Newton polytope of
$\Delta_A$ is a Minkowski summand of the secondary polytope of $A$,
which in turn implies that the secondary fan $\Sigma(A)$ is a
refinement of the normal fan of the Newton polytope of $\Delta_A$. The
tropical descriminant $\tau(X_A^*)$ being the codimension 1-skeleton of the
normal fan of the Newton polytope of $\Delta_A$, we obtain that
$\tau(X_A^*)$ is a subfan of the secondary fan $\Sigma(A)$. 

For a given $w\,{\in}\,\tau(X_A^*)$, the corresponding regular
subdivision $\Pi_w$, i.e., the support cell in $\Sigma(A)$, is
obtained as follows: Lift the point configuration $a_1,\ldots,a_n$
from $\R^d$ into $\R^{d+1}$ by extending with the coordinates of $w$
to $(a_1,w_1), \ldots, (a_n,w_n)$ in $\R^{d+1}$. The cells of $\Pi_w$
are the subsets of $A$ corresponding to the lower facets of the convex
hull of the lifted point configuration.
   
In fact, we conjecture that membership in $\tau(X_A^*)$ depends only
on the regular subdivision specified by the vector $w$, even in the
defective case:

\begin{conjecture} \label{conj:secfan}
For any point configuration $A$, the tropical discriminant 
$\tau(X_A^*)$ is a union of cones in the secondary fan~$\Sigma(A)$.
\end{conjecture}

Another notion that arose in the work of Gel'fand, Kapranov and Zelevinsky 
is the notion of {\em $\Delta$-equivalence\/} of regular triangulations of a 
non-defective point configuration: 
Let $\Pi_w$ and $\Pi_{w'}$ be two regular triangulations
which are neighbors in the secondary fan $\Sigma(A)$. This means that
their cones in $\Sigma(A)$ share a common face of codimension one. We
call $\Pi_w$ and $\Pi_{w'}$ {\em $\Delta$-equivalent\/} if they
specify the same leading monomial of the $A$-discriminant, i.e.,
$\,{\rm in}_w(\Delta_A) \,= \,{\rm in}_{w'}(\Delta_A)$.
The $\Delta$-equivalence classes of regular triangulations
of a point configuration $A$ define  a partition of the set of
maximal cones in the secondary fan $\Sigma(A)$. 

\begin{remark}
The $\Delta$-equivalence classes are in bijection with the
connected components of the complement  $\,\R^n \backslash \tau(X_A^*)\,$
of the tropical discriminant. 
\end{remark}

We return to Example~\ref{ex:symm3x3} to illustrate the relation of the 
tropical discriminant and the secondary fan of a point configuration as well 
as the notion of $\Delta$-equivalence.

\begin{figure}[ht]
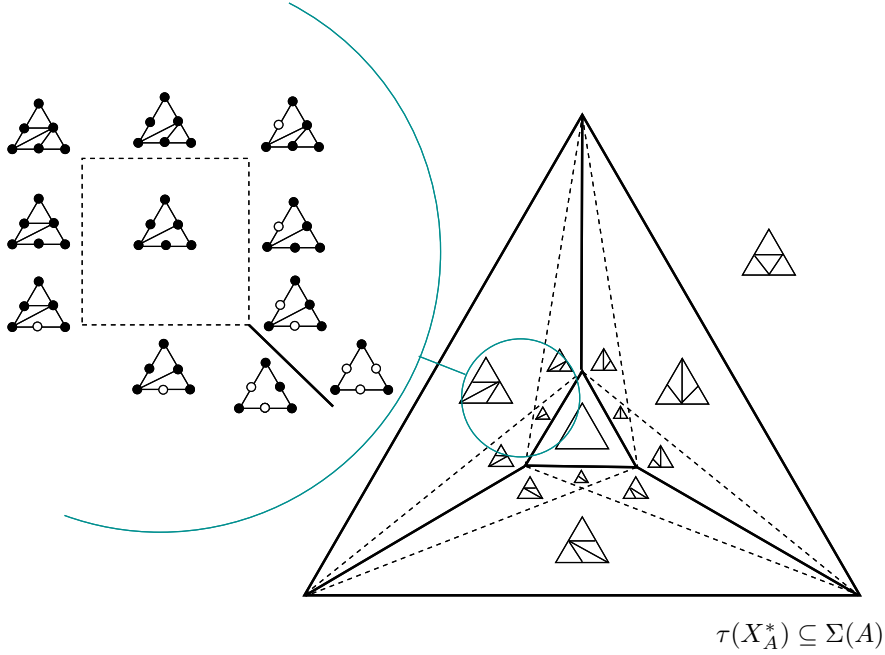

   \begin{picture}(0,0)%
     \includegraphics{secfan.pstex}%
   \end{picture}%
   \input{secfan.pstex_t}%
   
\caption{The tropical discriminant is
a subfan of the secondary fan} 
\label{fig_secfan3x3}
\end{figure}

\begin{example} \label{ex:symm3x3rev}
Let $A$ be the $3 \times 6$-matrix in
Example \ref{ex:symm3x3} whose toric variety
is the Veronese surface in $\C \P^5$ and
whose tropical discriminant $\tau(X_A^*)$
was depicted in Figure \ref{fig_B3x3}.
We now identify $\tau(X_A^*)$ as a subfan of
the secondary fan $\Sigma(A)$. Both of these two-dimensional fans
are drawn as planar graphs in Figure~\ref{fig_secfan3x3}. 
The graph $\Sigma(A)$ is dual to the three-dimensional
associahedron, and its $14$ regions are labeled
with the $14$ regular triangulations of the configuration~$A$.
The tropical discriminant $\tau(X_A^*)$ is the subgraph which is indicated
by solid lines. Edges of $\Sigma(A)$ that 
do not belong to $\tau(X_A^*)$ are dashed.
The  $14$ regular triangulations of~$A$ 
occur in five $\Delta$-equivalence classes
corresponding to the open cells in the complement of $\tau(X_A^*)$.

The magnification on the left in Figure~\ref{fig_secfan3x3} shows a portion
of the (dual) secondary polytope of $A$ with corresponding polyhedral 
subdivisions of the 6-point configuration indicated next to the faces.
\end{example}

\smallskip

We present two more examples which illustrate the 
various natural fan structures on the tropical discriminant. In
the non-defective case, the fan structure inherited from the secondary
fan refines the fan structure dual to the Newton polytope of $\Delta_A$.
The following example shows that a proper subdivision can occur.

\begin{example} \label{ex:K_4}
 Let  $d=3, n=6$, and consider the non-defective configuration
$$ A \quad = \quad\begin{pmatrix}
              1 & 1 & 1 & 1 & 0 & 0 \\
              0 & 0 & 0 & 0 & 1 & 1 \\
              0 & 1 & 2 & 3 & 0 & 1
\end{pmatrix}. $$
Here, $\,\Delta_A = x_1 x_6^3  - x_2 x_5 x_6^2 + x_3 x_5^2 x_6 - x_4 x_5^3 $.
Modulo the row space of $A$, the codimension 1-skeleton of the
secondary fan is two-dimensional. It is the
 the cone over a planar graph with
eight vertices and $18$ edges. The tropical discriminant
$\tau(X_A^*)$ corresponds to the induced subgraph on six of the
vertices, namely, the images of the six dimensional weights
$e_1, e_2, e_3, e_4, 2 e_1+e_2$ and $e_3+2 e_4$.
Here, {\em the secondary fan strictly refines the Gr{\"o}bner fan
on $\tau(X_A^*)$}. The latter is combinatorially a complete graph $K_4$ with
vertices $e_1, e_2, e_3 $ and $e_4$, while the former has the edges
${e_1,e_2}$ and ${e_3,e_4}$ subdivided by the vertices $2 e_1 + e_2$
and $e_3 + 2e_4$, respectively.
\end{example}

For defective configurations, the Gr\"obner fan structure can properly 
refine the secondary fan structure on 
$\tau(X_A^*)$, as the following example shows.

\begin{example}\label{ex:BJSST}
 We take $d=4, n=9$, and consider the defective configuration
$$ A \quad = \quad 
\begin{pmatrix}
              1 &1 & 1 & 0 & 0 & 0 & 0 & 0 & 0 \\ 0 & 0 &0& 1 &1 & 1 &
              0 & 0 & 0 \\ 0& 0 & 0 & 0 & 0 &0 & 1 & 1 & 1 \\ 0 & 1 &
              2& 0 & 1 & 2 & 0 & 1 & 2 \end{pmatrix}.  $$ 
The tropical discriminant $\tau(X_A^*)$ is a $7$-dimensional fan in
$\R^9$, regarded as a $6$-dimensional polyhedral complex.
Combinatorially, it is an immersion of the complete tripartite
hypergraph $K_{3,3,3}$. {\em The Gr{\"o}bner fan subdivision} has $51$
maximal cones and it {\em strictly refines the secondary fan
subdivision} which has only $49$ cones. Indeed, the vector $w = (0 ,1
,0 ,0 ,1 ,0 ,0 ,1 ,0)$ lies in the relative interior of a maximal cone
of the secondary fan subdivision which breaks into three maximal cones
in the Gr{\"o}bner fan subdivision. This example was verified by
applying the software {\tt Gfan} \cite{BJSST} to the
equations defining $X_A^*$.  In particular, we found that ${\rm
in}_w(X_A^*)$ is the codimension two primary ideal generated by the
determinant of the $3 \times 3$-matrix $$ \begin{pmatrix} x_1 & x_4 &
x_7 \\ x_2 & x_5 & x_8 \\ x_3 & x_6 & x_9
\end{pmatrix}      $$                                                         
plus the square of the ideal of $2 \times 2$-minors of the $2 \times3$-matrix  
$\, \begin{pmatrix}
x_1 & x_4  & x_7 \\
x_3 & x_6 & x_9
\end{pmatrix}$.
\end{example}

\medskip

There is yet another way to interpret the $A$-discriminant and its
tropicalization. For $A\in\Z^{d\times n},$ under the conditions imposed
above, the $A$-discriminant is the closure of the set of points
$[x_1: \cdots:x_n]$ in $\C\P^{n-1}$ such that the hypersurface in $(\C^*)^d$ 
defined by the Laurent polynomial
\[
      f(t)\, = \, \sum_{i=1}^n\, x_i\, t^{a_i}
\]
has a singular point in $(\C^*)^d$.

 In the tropical world, for $w\in
\R^n$, the regular polyhedral subdivision $\Pi_w$ of $A$ is
geometrically dual to the tropical hypersurface defined by the
tropical polynomial $\,\bigoplus_{i=1}^n w_i \odot t^{\odot a_i}
$. Hence, we may consider such tropical hypersurface to be {\em
singular\/} whenever the vector $w$ lies in the tropical discriminant
$\tau(X_A^*)$.  In this sense, our object of study in this article is
the {\em space of singular tropical hypersurfaces\/} and Theorem
\ref{thm:initial} gives a refined formula for the degree of that
space.

\begin{remark}
The tropical discriminant $\tau(X_A^*)$ is the polyhedral space
which pa\-ram\-e\-trizes
 all singular tropical hypersurfaces with fixed monomial support $A$.
\end{remark}

This relates our results to the celebrated work of Mikhalkin \cite{M}
on Gromov-Witten invariants.  For $d\,{=}\,3$, the $A$-discriminant
$X_A^*$ parametrizes singular curves on the toric surface~$X_A$.  Our
formula for the degree of $X_A^*$ is consistent with the lattice paths
count in \cite{M} for the number of nodal curves on $X_A$ of genus
$g{-}1$ through $n{-}2$ fixed points, where $g$ is the number of
interior points in the lattice polytope~$Q_A$.  It would be
interesting to explore possible applications of our combinatorial
approach to Gromov-Witten theory.  The work of Gathmann and Markwig
\cite{GM} offers an algebraic setting for such a study (see also the
recent preprint \cite{EK}, which addressed this question after the first
version of our paper).


\section{Cayley configurations and resultant varieties}
\label{sec:resultants}

One of the main applications of $A$-discriminants
is the study of resultants in {\em elimination theory}.
The configurations $A$ which arise in elimination theory
have a special combinatorial structure arising from 
the {\em Cayley trick}. See
 \cite[\S 3.2.D]{GKZ} for a geometric introduction.
Based on the results of the earlier sections,
we here study the combinatorics and geometry of
{\em tropical resultants}, and we generalize the positive
degree formula for resultants 
in \cite{BS} to
resultant varieties of arbitrary codimension.

Let $A_1,\dots,A_m$ be finite subsets of $\Z^r$. Their {\em
   Cayley configuration} is defined as
\begin{equation}
\label{Cayley}
   A \,\,\, = \,\, \,
  \{e_1\} \! \times A_1 \! \, \, \cup \,\,
  \{e_2\} \! \times A_2\! \, \, \cup \,\,
\cdots   \, \, \cup \,\,
\{e_{m}\}   \! \times \! A_{m}
\,\, \subset \,\, \Z^{m} \times \Z^r ,
\end{equation}
where $e_1,\dots, e_m$ is the standard basis of $\Z^{m}$.
To be consistent with our notation in Sections 1--5, we can regard
$A$ as a $d \times n$-matrix with
$d = m+r$ and $\,n = |A_1| + |A_2| + \cdots + |A_m|$.
As in \cite[\S 8.1]{GKZ} and \cite{BS},  the Cayley configuration $A$ 
represents the following system 
of $m$ Laurent polynomial equations in $r$ unknowns:
\begin{equation}
\label{Cayleysystem}
\sum_{u \in A_1} x_{1,u} z^u  \,\, = \,\,
\sum_{u \in A_2} x_{2,u} z^u  \,\, = \,\,
\cdots  \,\, = \,\,
\sum_{u \in A_m} x_{m,u} z^u  \,\,\, = \,\,\, 0 .
\end{equation}
Here $z = (z_1,\ldots,z_r)$ are coordinates
on $(\C^*)^r$ and we use multi-index notation
$\,z^u = z_1^{u_1}  \cdots z_r^{u_r}$.
Our earlier examples include the following Cayley configurations:

\begin{itemize}
\item In Example  \ref{mixeddiscriminant} we have
$m = 2, r = 2$, and the system (\ref{Cayleysystem}) takes the form
\begin{eqnarray*}
& \,\, x_1 z_1^2 z_2^{19}  \,+\, 
x_2 z_1^3 z_2^{17} \,+\, 
x_3 z_1^5 z_2^{13}   \,+\, 
x_4 z_1^7 z_2^{11}  \,\,\,
\quad = \quad 0 , \\ &
x_5  z_1^{11} z_2^7  + 
x_6  z_1^{13} z_2^5  + 
x_7 z_1^{17} z_2^3   + 
x_8 z_1^{19} z_2^2  
\quad = \quad 0. 
\end{eqnarray*}
\item In  Example \ref{ex:K_4} we have 
$m=2$, $r=1$, $A_1 = \{0,1,2,3\}$ and $A_2 = \{ 0,1 \}$.
The $A$-discriminant is the Sylvester resultant
\[
\Delta_A ={\rm det}
\begin{pmatrix}
x_1 & x_2 & x_3 & x_4 \\
x_5 & x_6 & 0   & 0  \\
 0  & x_5 & x_6 & 0   \\
 0  &  0  & x_5 & x_6  
\end{pmatrix}\, .
\]
\item In  Example \ref{ex:BJSST}
we have $m=3$, $r= 1$ and
$A_1 = A_2 = A_3 = \{0,1,2\}$, and
the system (\ref{Cayleysystem}) consists of
three quadratic equations in one unknown $z$:
$$ x_1 + x_2 z + x_3 z^2 \,\, = \,\,
 x_4 + x_5 z + x_6 z^2 \,\, = \,\,
 x_7 + x_8 z + x_9 z^2 \,\,\, = \,\,\, 0 .$$
 The variety $X_A^*$  of all solvable systems
 of this form has codimension two.
 \end{itemize}

\smallskip

Returning to the general case, we say that the
Cayley configuration $A$ is {\em essential} 
if the Minkowski sum
$\,\sum_{i \in I} A_i \,$ has affine dimension at 
least $|I|$ for every subset~$I$ of $\{1,\dots,m\}$ with $|I|\leq r$.
The {\em resultant variety}
of the Cayley configuration $A$ is the Zariski closure
in $\C \P^{n-1}$ of the set of all points $(x_1:x_2:\ldots:x_n)$
whose corresponding system (\ref{Cayleysystem}) has a solution
$z$ in $(\C^*)^r$.
The following result is a generalization of 
Proposition 1.7 in \cite[\S 9.1.A]{GKZ} and of
Proposition 5.1 in \cite{CDS}. 

\begin{proposition} \label{prop:cayley}
The resultant variety of any Cayley configuration $A$
contains the $A$-discriminant $X_A^*$.
If $m \geq r+1$ and the Cayley configuration $A$ is essential then
the resultant variety and the $A$-discriminant
coincide.
\end{proposition}

\begin{proof}
Consider the hypersurface in $(\C^*)^m \times (\C^*)^r$
defined by the equation
$$ \sum_{i=1}^m \sum_{u \in A_i} x_{i,u} \cdot t_i z^u  \,\, = \,\, 0 . $$
If $(t,z) \in (\C^*)^{m+r}$ is a singular point on this
hypersurface then $z \in (\C^*)^r$ is a solution to
(\ref{Cayleysystem}). This proves the inclusion.
If $A$ is essential then a linear algebra argument
as in \cite[\S 9.1.A]{GKZ} shows that every
solution $z$ of (\ref{Cayleysystem}) arises in this way.
\end{proof}

The hypothesis that $A$ be essential is necessary
for the equality of the resultant variety and the
$A$-discriminant, even when $m=r+1$, the situation
of  classical elimination theory. The following simple 
example illustrates the general behavior.

\begin{example}\label{ex:nonessential}
Let  $r=2, m=3$, $A_1 = A_2 = \{ (0,0), (1,0)\}$,
and $A_3 = \{ (0,0), (1,0),$ $(0,1),(1,1)\}$. The
Cayley matrix $A$ represents a toric
$4$-fold $X_A$ in $\C \P^7$. 
It is not essential since
$ A_1  +  A_2  $ is one-dimensional.
The system (\ref{Cayleysystem}) equals
$$  x_1 + x_2 z_1 \,= \, x_3 + x_4 z_1 \,= \,
x_5 + x_6 z_1 + x_7 z_2 + x_8 z_1 z_2 \,\,\, = \,\,\, 0 . $$
The resultant variety has codimension one, with equation 
$\,x_1 x_4 = x_2 x_3$, but the $A$-discriminant $X_A^*$
has codimension three. In fact, we have $X_A^* = X_A$
in this case.
\end{example}

For the rest of this section we consider an essential
Cayley configuration as in (\ref{Cayley}) with  $m \geq r+1$ blocks and we
set $c =m-r$. Then the following result holds.

\begin{lemma} \label{rescodimc}
The resultant  variety $X_A^*$ has codimension $c$.
\end{lemma}

\begin{proof}
Let $W$ denote the incidence variety consisting
of all pairs $(x,z)$ in
$\,\C\P^{n-1}  \times (\C^*)^r\,$ such that
(\ref{Cayleysystem}) holds.
Let $\pi_1: W \to \C\P^{n-1}$  be the  projection to the first factor.
By Proposition~\ref{prop:cayley},
the resultant variety $X_A^*$ coincides with the closure of $\pi_1(W)$.
Looking at the second projection $\pi_2: W \to (\C^*)^r$,
which is surjective and whose fibers are linear spaces of dimension $n-1 -m$,
we deduce that
$W$ is irreducible and  has dimension  $(n-1-m)+r = n-1-c$. 
Then, $\dim(X_A^*) \leq n-1-c$.

Given a generic choice of coefficients $x_i$
for the first $r$ polynomials,
it follows from the essential hypothesis  and Bernstein's Theorem,
that the first $r$ equations in (\ref{Cayleysystem}) have
a common solution $z \in (\C^*)^r$. We can freely
choose all but one of the coefficients of the last $c$ polynomials so that $z$
solves  (\ref{Cayleysystem}). This implies that
 $\dim(X_A^*) \geq n-1-c$, and the lemma follows.
\end{proof}

Corollary~\ref{cor:dimension} asserts that there exists a 
chain $\sigma_1, \dots, \sigma_{n- 2m}$ of $(0,1)$-vectors
representing the supports of vectors in the kernel of $A$ such that
the rank of the matrix $ (A^t \!,\sigma_1, \dots, \sigma_{n-2m})$ is
precisely $n-c$.  We present an explicit way of choosing such a
chain. By performing row operations, we can assume that each
set $A_i$ contains the origin. Set $B_i = A_i \backslash \{0\}$. Let
$b_i \in B_i$ for $i=1,\dots, r$ such that $b_1, \dots, b_r$ are
linearly independent. Such elements exist because the family of
supports is essential. Now, for any other element $a$ in $B = (B_1
\cup \dots \cup B_{r+c}) \backslash \{ b_1, \dots, b_r\}$, we can find
an element $v_a$ in $\ker(A)$ with support corresponding to the origin
in each $A_i$ for $i$ from~$1$ to $r$, union the variables
corresponding to $b_1, \dots, b_r$ and $a \in B_j$, plus the origin of
$A_j$ in case $j >r$. Choose any such $a \in B_{r+1}$ and let
$\sigma_1$ be the support of $v_a$; it will have $2r +2$ non-zero
coordinates. Add a new point $a'$ in $B$. We can assume that the
support of $v_a + v_a'$ equals the union of their supports. 
Let~$\sigma_2$ be the associated support vector. We continue in this
manner, adding a new point in $B$ at a time, and considering a new
element in the chain of support vectors, but avoiding to pick all of
$B_1 \cup \dots \cup B_{r+1}$ and all of each of $B_{r+2}, \dots,
B_{r+c}$.  This produces precisely
$$1 + (n_1 -2 + \dots + n_{r+1} -2)-1 + (n_{r+2} -2) + \dots +
(n_{r+c}-2) = n -2r - 2 c = n -2m$$
vectors $\sigma_1, \dots,
\sigma_{n-2m}$ in ${\mathcal C}(A)$.  It is straightforward to check
that the rank of the submatrix of $ (A^t \!,\sigma_1, \dots,
\sigma_{n-2m})$ given by the first $m$ and the last $n -2m$ columns
has maximal rank $n-m$.  Note that this is just a $(0,1)$ matrix.
Adding the last $r$ columns of $A^t$ containing the information about
the specific supports $A_1, \dots, A_m$ will increase the rank by $r$,
as a consequence of the fact that the family is essential. Therefore,
the rank of the matrix $ (A^t \!,\sigma_1, \dots, \sigma_{n-2m})$ is
precisely $n - m + r = n-c$.

\smallskip

We identify $\{1,2,\ldots,n\}$ with the disjoint
union of the sets $A_1,A_2,\ldots,A_m$. Thus
a generic vector $w \in \R^n$ assigns a height
to each point in any of the $A_i$, and it
defines a {\em tight coherent mixed subdivision} $\Delta_w$
of the Minkowski sum $\sum_{i=1}^m A_i$ (cf. \cite{BS}).
When $c=1$,  the initial form with respect to $w$ of the 
mixed resultant
is described in \cite[Theorem 2.1]{BS} in terms of the sum of volumes
of suitable {\em mixed cells} of  the tight
coherent mixed decomposition (TCMD) induced by
$w$. We next generalize this result to 
resultant varieties of arbitrary codimension $c$. 
For a classical study of resultant ideals  of
dense homogeneous polynomials we refer to \cite{Jou}.

\begin{theorem} \label{thm:Sec6}
A codimension $c$ monomial prime $P_\tau = \langle x_{\tau_1}, \dots,
x_{\tau_c} \rangle$
is a minimal prime of the monomial ideal
$\,{\rm in}_w(X_A^*)\,$ only if $\tau$ consists
of one point each from $c$ of the $A_j$. The
multiplicity of $P_\tau$ is the  total volume of
all mixed cells  in
the tight coherent mixed subdivision $\,\Delta_w \,$ which use the
points of $\tau$ in their decomposition.
\end{theorem}

\begin{proof}
The resultant variety $X_A^*$ is irreducible
(by Proposition \ref{prop:cayley}), and it has  codimension $c$
(by Lemma \ref{rescodimc}). This implies 
(by \cite[Theorem 1]{KS}) that every minimal prime of
the initial monomial ideal 
$\,{\rm in}_w(X_A^*)\,$ has codimension $c$.
Let $\,P_\tau = \langle x_{\tau_1}, \dots,
x_{\tau_c} \rangle\,$ be such a minimal prime.
After relabeling we may assume that
each $x_{\tau_i}$ in $P_\tau$ is a coefficient of
one of the last $c$ Laurent polynomials in
(\ref{Cayleysystem}).

The proof for the case $c=1$ is given in
\cite[\S 2]{BS}, and the proof for $c > 1$
uses the same general technique.
We write $f_i$ for the $i$-th Laurent polynomial in (\ref{Cayleysystem}),
but with $x_i$ replaced by $x_i \epsilon^{w_i}$.
Let ${\bf K}$ be the algebraic closure of the field of rational
functions over $\,\C\{\!\{\epsilon\}\!\}\,$
 in the coefficients of the first
$r$ Laurent polynomials $f_1,\ldots,f_r$, 
let ${\bf x}$ denote the vector all coefficients
of the last $c$ Laurent polynomials $f_{r+1},\ldots,f_m$, 
and consider the polynomial ring ${\bf K}[{\bf x}]$.
Let $\mu$ denote the mixed volume of the Newton polytopes
of the polynomials $f_1,\ldots,f_r$. Then, by
Bernstein's Theorem, the system
$f_1 = \cdots = f_r=0$ has $\mu$ distinct roots
${\bf z}_1(\epsilon),\ldots,{\bf z}_\mu(\epsilon)$ in $({\bf K}^*)^r$.

For any $j \in \{1,2,\ldots,\mu\}$, the ideal
$\,I_j \,=\,  \langle f_{r+1}({\bf z}_i(\epsilon)),\ldots,
f_{m}({\bf z}_i(\epsilon))\rangle \,$
is generated by linear forms in ${\bf K}[{\bf x}]$.
The intersection of these ideals,
$\,I\, =  \,I_1 \cap I_2 \cap \cdots \cap I_\mu$,
is an ideal of codimension $c$ and degree $\mu$
in ${\bf K}[{\bf x}]$.
Geometrically, we obtain $I$ by embedding the prime ideal of
 $X_A^*$ into ${\bf K}[{\bf x}]$ and then
replacing $x_i$ by $x_i \epsilon^{w_i}$.
This is the higher codimension version of the
{\em product formula for resultants} \cite[Eqn.~(14)]{BS}.

The ideal $I$ represents a flat family,
and its special fiber $I|_{\epsilon=0}$ 
at $\epsilon = 0$ coincides
with the special fiber of  the image of ${\rm in}_w(X_A^*)$ 
in ${\bf K}[{\bf x}]$. In particular, $P_\tau$
is an associated prime of $I|_{\epsilon=0}$,
and it contains one of the
ideals $\,I_j|_{\epsilon=0}$. Since
the generators of $I_j$ are $c$ linear forms
in  disjoint groups of unknowns $x_\ell$,
we see that $P_\tau$ contains one
unknown from each group. This proves
the first statement in Theorem \ref{thm:Sec6}.

 After relabeling
we may assume that $x_{\tau_j}$ is
a coefficient of $f_{r+j}$ for $j=1,\ldots,c$.
Each root ${\bf z}_j(\epsilon)$ corresponds
to a mixed cell $C$ in the TCMD of the 
small Minkowski sum $A_1+\cdots+A_r$ defined by the
restriction of $w$.
By the genericity of $w$,
the mixed cell $C$ corresponds to a unique
cell $C'$ in the TCMD $\Delta_w$ of the
big Minkowski sum $A_1+\cdots+A_r+A_{r+1} +\cdots + A_m$,
and every mixed cell of $\Delta_w$ arises in this manner.
The reasoning above implies that
the mixed cell $C'$ uses the points of $\tau$ in its decomposition
if and only if  $\,I_j|_{\epsilon=0} = P_\tau\,$ in
${\bf K}[{\bf x}]$. This completes the proof.
\end{proof}

The first assertion in Theorem  \ref{thm:Sec6}
can also be derived more easily, namely, from the fact that
for any $(r+1)$-element subset $I$ of  $\,\{ 1, \ldots, m\}$,
the mixed resultant of the configurations $A_{i}, i \in I,$
vanishes on $X_A^*$ and only contains unknowns
 $x_{i,a}$ with $i \in I$. However, for
the multiplicity count in the second assertion
we need the ``product formula'' developed above.
Theorem \ref{thm:Sec6} has the following corollary.

\begin{corollary} \label{cor:Sec6}
The degree of the resultant variety $X_A^*$ is the
sum of the mixed volumes $MV(A_{i_1}, \ldots,A_{i_r})$
as $\{i_1,\ldots,i_r\}$ runs over all
$r$-element subsets of $\{1,\ldots,n\}$.
\end{corollary}

We present two examples to illustrate Theorem \ref{thm:Sec6}
and Corollary \ref{cor:Sec6}.

\begin{example}
Let $m=3$, $r= 1$ and $A_1 = A_2 = A_3 = \{0,1,2\}$
as in Example \ref{ex:BJSST}, and choose
$w \in \R^9$ which represents the
{\em reverse lexicographic term order}. Then 
$$ {\rm in}_w(X_A^*) \, =  \, \langle \,
x_3 x_5 x_7 , \,
x_6^2 x_7^2 , \,
x_3 x_6 x_7^2 , \,
x_3^2 x_7^2 , \,
x_3 x_4 x_6 x_7 , \,
x_3^2 x_4 x_7 , \,
x_3^2 x_4^2 , \,
x_2 x_4 x_6^2 x_7\,\rangle . $$
This ideal has seven associated primes,
of which three are minimal:
$\,\langle x_3, x_6 \rangle$, $\,\langle x_3, x_7 \rangle$, 
and $\,\langle x_4, x_7 \rangle$.
They correspond
to the  three mixed cells  
$\,(3,6, \{7,8,9\})$,
$\,(3,\{4,5,6\},7)\,$
and $\,(\{1,2,3\}, 4,7)\,$ of
the TCMD $\Delta_w$ of
$\,A_1 + A_2 + A_3 = \{0,1,\ldots,6\}$.
Each mixed cell has volume two, which 
implies that the degree of $X_A^*$ is $2+2+2 = 6$.
\end{example}

\begin{example}
Let $r=2$ and $m=4$ and take the $A_i$ to be
the four subtriangles of the square
with vertices $(0,0)$, $(0,1)$, $(1,0)$ and $(1,1)$.
Here  (\ref{Cayleysystem}) is a system of four equations in two
unknowns $z_1$ and $z_2$ which can be written
in matrix form as 
\begin{equation} \label{eq:4by4}
\begin{pmatrix}
x_{1}  & x_{2} & x_{3} &  0   \\
x_{4}  & x_{5} & 0 & x_{6} \\
x_{7}  &  0  & x_{8} & x_{9} \\
  0 & x_{10} & x_{11} & x_{12} \\
\end{pmatrix}
\cdot
\begin{pmatrix}
1 \\ z_1 \\ z_2 \\ z_1 z_2
\end{pmatrix}
\quad = \quad
\begin{pmatrix}
0 \\ 0 \\ 0 \\  0 \end{pmatrix}.
\end{equation}
The resultant variety $X_A^* \subset \C\P^{11}$ has codimension $2$ 
and degree  $\, 12 = \binom{4}{2} \cdot 2$.
This is the sum over the mixed areas
of the $\binom{4}{2}$ pairs of triangles.
The mixed area is~$2$ for each such pair.
Using computer algebra, we find that 
the prime ideal of the resultant variety $X_A^*$ is generated
by the $4 \times 4$ determinant in (\ref{eq:4by4})
together with ten additional polynomials of degree six
in the $x_i$. Theorem \ref{thm:Sec6}
gives a combinatorial recipe for constructing all
the initial monomial ideals of this prime ideal.
\end{example}

We conclude with a brief discussion
of the {\em tropical resultant} $\tau(X_A^*)$.
The results in Section~\ref{sec:subdivisions} characterize this polyhedral
fan in terms of regular subdivisions of $A$, and we
will now rephrase this characterization in terms of
coherent mixed subdivisions of $(A_1,\ldots,A_m)$.
Theorem 5.1 in \cite{BS} implies  that every
regular subdivision $\Pi_w$ of~$A$
corresponds uniquely to a coherent mixed subdivision
(CMD), which we denote by $\Delta_w$. Note that $\Pi_w$ is a polyhedral
cell complex of dimension $r+m-1$ while
$\Delta_w$ has only dimension $r$, and
$\Pi_w$ is a triangulation if and only if $\Delta_w$ is a TCMD.

Every cell $F$ of a CMD $\Delta_w$ decomposes uniquely
as a Minkowski sum $\, F = F_1 + \cdots + F_m$,
where $F_i \subset A_i$ for all $i$. We write $\overline{F}$
for the corresponding cell of $\Delta_w$.
We say that the cell $F$ is {\em fully mixed} 
if each $F_i$ has affine dimension at least one.
Using the techniques in the proof of Theorem \ref{thm:Sec6}, 
we can derive the following:

\begin{proposition} \label{prop:subd}
The tropical resultant $\,\tau(X_A^*)  \,$ 
equals the set
\begin{equation}
\label{IsFullyMixed}
\bigl\{ w \in \R^n \,:\,
\hbox{ $\Delta_w$ has a maximal cell $F$ which is fully mixed }
\bigr\}. 
\end{equation}
\end{proposition}

Specializing to the classical case $c=1$, when
$X_A^*$ is the hypersurface defined by the
mixed resultant, we can now recover the
combinatorial results in \cite[\S 5]{BS}. In particular,
the positive formula for the extreme monomials of the
sparse resultant given in \cite[Theorem 2.1]{BS}
can now be recovered as a special case of Theorem 
\ref{thm:initial}. Thus, one way to look at Sections 1--5
in the present paper is that these extend all  results
in \cite{BS} from essential Cayley configurations 
with $m=r+1$ to arbitrary configurations $A$.
A perspective on how this relates to the approach of
 \cite{GKZ} is given by the
points (a),(b),(c),(d) found at the end
of the introduction in \cite[p.~208-209]{BS}.


\bibliographystyle{amsplain}

\end{document}